\documentclass[preprint,11pt]{elsarticle}
\usepackage[toc,page]{appendix}

\usepackage{amsmath,amsthm,amsfonts,mathtools,bm,hyperref}
\usepackage{amssymb}
\usepackage{dsfont}
\usepackage{mathrsfs}
\usepackage{multirow}
\usepackage{tikz,pgfplots}
\usetikzlibrary{patterns}
\usepackage{subcaption}
\pgfplotsset{compat=newest}
\usepackage{float,indentfirst,algorithm}
\usepackage{pgfplots}
\usepackage{booktabs}
\usepackage{hyperref}
\usepackage{multicol}
\usepackage{nomencl}
\usepackage[margin=1in]{geometry}
\usetikzlibrary{patterns}

\makenomenclature

\usepackage{etoolbox}
\renewcommand\nomgroup[1]{%
  \item[\bfseries
  \ifstrequal{#1}{A}{Symbols}{%
  \ifstrequal{#1}{B}{Operators}{}}%
]}

\newtheorem{theorem}{Theorem}
\newtheorem{lemma}[theorem]{Lemma}
\theoremstyle{plain}

\newtheorem{remark}{Remark}
\newtheorem{definition}{Definition}
\usepackage[capitalise]{cleveref}
\usepackage{appendix}
\usepackage{multicol}
\usepackage{nomencl}

\newcommand{\mean}{m}
\newcommand{\pmean}{\widehat{m}}
\newcommand{\pp}{\widehat{\theta}}
\newcommand{\py}{\widehat{y}}
\newcommand{\ppy}{\widehat{y}}
\newcommand{\Cov}{C}
\newcommand{\pCov}{\widehat{C}}
\newcommand{\CovZ}{Z}
\newcommand{\pCovZ}{\widehat{Z}}

\newcommand{\pGCovZ}{\widehat{\mathcal{Y}}}

\newcommand{\N}{\mathcal{N}}
\newcommand{\G}{\mathcal{G}}
\newcommand{\F}{\mathcal{F}}
\newcommand{\E}{\mathbb{E}}
\newcommand{\I}{\mathbb{I}}
\newcommand{\R}{\mathbb{R}}

\newcommand{\bigO}{\mathcal{O}}

\makeatletter
\def\ps@pprintTitle{%
   \let\@oddhead\@empty
   \let\@evenhead\@empty
   \def\@oddfoot{\reset@font\hfil\thepage\hfil}
   \let\@evenfoot\@oddfoot
}
\makeatother

\begin{document}

\begin{abstract}
The unscented Kalman inversion~(UKI) presented in~\cite{UKI} is a general derivative-free approach to solving the inverse problem. 
UKI is particularly suitable for inverse problems where the forward model is given as a black box and may not be differentiable.
The regularization strategy and convergence property of the UKI are thoroughly studied, and the method is demonstrated effectively handling noisy observation data and solving chaotic inverse problems. 
In this paper, we aim to make the UKI more efficient in terms of computational and memory costs for large scale inverse problems. We take advantages of the low-rank covariance structure to reduce the number of forward problem evaluations and the memory cost, related to the need to propagate large covariance matrices. And we leverage reduced-order model techniques to further speed up these forward evaluations.
The effectiveness of the enhanced UKI is demonstrated on a barotropic model inverse problem with $\bigO(10^5)$ unknown parameters and a 3D generalized circulation model~(GCM) inverse problem, where each iteration is as efficient as that of gradient-based optimization methods.
\end{abstract}

\begin{keyword}
  Inverse Problem, Optimization, Kalman Filter, Low-Rank Approximation, Reduced-Order Model
\end{keyword}

\begin{frontmatter}

  \title{Improve Unscented Kalman Inversion With Low-Rank Approximation and Reduced-Order Model}

  \author[rvt1]{Daniel~Z.~Huang}
  \ead{dzhuang@caltech.edu}
  
  \author[rvt2]{Jiaoyang~Huang}
  \ead{jh4427@nyu.edu}

  \address[rvt1]{California Institute of Technology, Pasadena, CA}
  \address[rvt2]{New York University, New York, NY}

\end{frontmatter}

\nomenclature[A]{$\theta$}{unknown parameter vector} 
\nomenclature[A]{$y$}{observation vector} 
\nomenclature[A]{$\py$}{observation vector mean in prediction step} 
\nomenclature[A]{$Y$}{observation set} 
\nomenclature[A]{$\G$}{mapping between parameter space and observation space} 
\nomenclature[A]{$\F$}{filter} 
\nomenclature[A]{$t$}{time} 
\nomenclature[A]{$N_\theta$}{dimension of unknown parameter vector} 
\nomenclature[A]{$N_y$}{dimension of observation vector} 
\nomenclature[A]{$\eta$}{observation error} 
\nomenclature[A]{$\Sigma_{\eta}$}{observation error covariance}
\nomenclature[A]{$\omega$}{artificial evolution error} 
\nomenclature[A]{$\Sigma_{\omega}$}{artificial evolution error covariance}
\nomenclature[A]{$\nu$}{artificial observation error} 
\nomenclature[A]{$\Sigma_{\nu}$}{artificial observation error covariance}
\nomenclature[A]{$\mean$}{conditional mean} 
\nomenclature[A]{$\pmean$}{conditional mean in prediction step}
\nomenclature[A]{$\Cov$}{conditional covariance} 
\nomenclature[A]{$\pCov$}{conditional covariance in prediction step}
\nomenclature[A]{$\Phi$}{least-square function}
\nomenclature[B]{$c$}{unscented points weight}
\nomenclature[B]{$W$}{unscented points weight}
\nomenclature[B]{$\kappa$}{unscented points weight}
\nomenclature[B]{$\alpha$}{unscented points weight}
\nomenclature[B]{$\beta$}{unscented points weight}
\nomenclature[B]{$\Box_{n}$}{$n$-th time step}
\nomenclature[B]{$\Box^{j}$}{$j$-th ensemble particle}
\nomenclature[B]{$\Box_{(i)}$}{$i$-th component}


\section{Introduction}
Inverse problems are ubiquitous in engineering and scientific applications. These include, to name only a few, global climate model calibration~\cite{sen2013global,schneider2017earth,dunbar2020calibration}, material constitutive relation calibration~\cite{huang2020learning,xu2021learning,avery2020computationally}, seismic inversion in geophysics~\cite{russell1988introduction,bunks1995multiscale}, and medical tomography~\cite{goncharsky2017iterative,hakkarainen2019undersampled}.
The associated forward problems, which may feature multiple scales or include chaotic and turbulent phenomena, can be very expensive.  
Moreover, the observational data is often noisy, and the inverse problem possibly involves a large number of unknown
parameters to recover and may be ill-posed.

Inverse problems can be formulated as recovering unknown parameters~(or states) $\theta \in \R^{N_{\theta}}$ from the noisy observation $y \in \R^{N_y}$, as following
\begin{equation}
\label{eq:KI}
    y = \G(\theta) + \eta,
\end{equation}
where $\G$ denotes a mapping between the parameter space and the observation space, and $\eta \sim \N(0,\Sigma_{\eta})$ denotes the observational error.
Central to both the optimization and probabilistic approaches to inversion is the regularized objective function $\Phi_R(\theta)$ defined by 
\begin{subequations}
\label{eq:KI2}
\begin{align}
\Phi_R(\theta) &:= \Phi(\theta)+\frac{1}{2}\lVert\Sigma_{0}^{-\frac{1}{2}}(\theta - r_0) \rVert^2,\\
    \Phi(\theta) &:= \frac{1}{2}\lVert\Sigma_{\eta}^{-\frac{1}{2}}(y - \G(\theta)) \rVert^2,
\end{align}
\end{subequations}
where $\Sigma_{\eta} \succ 0 $ is strictly positive-definite
and normalizes the model-data misfit $\Phi$. 
$r_0$ and $\Sigma_{0}\succ 0$ encode
prior mean and covariance information about $\theta$.

The focus of this paper is mainly on the unscented Kalman inversion~\cite{wan2000unscented,UKI},  a version of the Kalman inversion method.
Kalman inversion methodology pairs the parameter-to-data map~$\G$ with a stochastic dynamical system for the parameter and then employs techniques from Kalman filtering~\cite{kalman1960new,sorenson1985kalman,evensen1994sequential,julier1995new} to estimate the parameter $\theta$ given the data $y$.
When the unscented Kalman filtering is applied, this leads to the UKI, a derivative-free method to solve the
inverse problem aimed at
solving the optimization problem   
defined by minimization of $\Phi_R$.
For inverse problems in general, UKI is attractive because it is derivative-free and hence introduces a significant flexibility in the forward solver design.
Therefore, UKI is suitable for complex multiphysics problems requiring coupling of different solvers~\cite{huang2019high,huang2020high,huang2018simulation,huang2020modeling,adcroft2019gfdl}, where forward models  can be treated  as  a  black  box,  and methods containing discontinuities~(i.e., immersed/embedded boundary method~\cite{peskin1977numerical,berger2012progress,huang2018family,huang2020embedded} and adaptive mesh refinement~\cite{berger1989local,borker2019mesh}). 
Moreover, 
Tikhonov type regularization embedded in the stochastic dynamics~\cite{UKI} enables UKI to effectively handle noisy observations; 
the Levenberg-Marquardt connection~\cite{UKI}  enables UKI  to achieve linear (sometimes superlinear) convergence  in practice;
and the smoothing property~\cite{UKI}  introduced by the unscented transform enables UKI to effectively handle chaotic inverse problems. 
Therefore, UKI  is a general approach to solving the inverse problem.

However, at  each  iteration,  the UKI  requires  evaluating  the  forward  problem  $2N_{\theta}+ 1$ times and propagating large $N_\theta$ by $N_\theta$ covariance matrices.   Although  these forward problem evaluations are embarrassingly parallel, the computational or memory costs can be intractable, when $N_{\theta}$ is large or the forward problem evaluation is expensive. The present work focuses on
alleviating the computational and memory costs of the UKI.

\subsection{Literature Review}
The Kalman filter~\cite{kalman1960new} and its variants, including but not limit to extended Kalman filter~\cite{sorenson1985kalman}, ensemble Kalman filter~\cite{evensen1994sequential,anderson2001ensemble,bishop2001adaptive}, unscented Kalman filter~\cite{julier1995new,wan2000unscented}, and cubature Kalman filter~\cite{arasaratnam2009cubature} are developed to sequentially update the probability distribution of states in  partially  observed  dynamics.
Most of Kalman filters use Gaussian ansatz to formulate Kalman gain to assimilate the observation and update the distribution.
The unscented Kalman filter represents the mean and covariance of the Gaussian distribution directly.
A quadrature rule---unscented transform, specifically for Gaussian distribution based on $2N_{\theta}+1$ deterministic $\sigma-$points~\cite{julier1995new,julier1997new,wan2000unscented}, has been designed and applied for the propagation of the state mean and covariance. 
Besides state estimation, Kalman  methods  are used for parameter inversion, which bring about the extended Kalman inversion~(ExKI), the ensemble Kalman inversion~(EKI), and the unscented Kalman inversion~(UKI). 
Extended, ensemble and unscented Kalman inversions have been applied to train neural networks~\cite{singhal1989training, puskorius1991decoupled, wan2000unscented, kovachki2019ensemble} and EKI has been applied in the oil industry~\cite{oliver2008inverse, chen2012ensemble, emerick2013investigation}. Dual and joint Kalman filters~\cite{wan1997neural, wan2000unscented} have been 
designed to simultaneously estimate the unknown states and the parameters~\cite{wan1997neural,parlos2001algorithmic,wan2000unscented,gove2006application,albers2017personalized} from noisy sequential observations.
The EKI has been systematically developed and analyzed as
 a general purpose methodology~\cite{iglesias2013ensemble,schillings2017analysis,schillings2018convergence,albers2019ensemble,chada2020tikhonov} for the solution of inverse and parameter
 estimation problems, the same is for UKI~\cite{UKI}.
However, unscented Kalman filter and inversion generally require $2N_{\theta}+1$ forward problem evaluations and the storage of a $N_{\theta} \times N_{\theta}$ covariance matrix, therefore, they are (mistakenly) supposed to be unsuitable for high-dimensional data-assimilation/inverse problems.

Since many physical phenomena or systems feature large-scale structure or low-dimensional attractors, therefore, 
the model error covariance matrices are generally low-rank.
The low-rank covariance structures are leveraged in ensemble based Kalman filters, which leads to reduced space Kalman filters~\cite{kaplan1997reduced,brasseur1999assimilation,gourdeau2000assimilation} and square root filters~(i.e, ensemble adjustment Kalman filter~\cite{anderson2001ensemble} and ensemble transform Kalman filter~\cite{bishop2001adaptive,tippett2003ensemble,wang2003comparison}). And these square root filters have been applied successfully for weather forecast.
Inspired by these square root filters~\cite{anderson2001ensemble,bishop2001adaptive,van2001square}, we introduce the reparameterization strategy and the low-rank square root form of the UKI, namely the truncated unscented Kalman inversion~(TUKI), where the inverse problems are solved in the vastly reduced subspace. This reduces the computation and memory cost from $\bigO(N_{\theta})$ and $\bigO(N_\theta^2)$ to $\bigO(N_r)$ and $\bigO(N_rN_{\theta})$, where $N_{r}$ is the rank of the covariance matrix.
It has also been pointed out in~\cite{anderson2001ensemble} that "filter methods are unlikely to work without making use of such information about the covariance structure".

Moreover, these $2N_{\theta}+1$ deterministic $\sigma-$points in the unscented transform representing the Gaussian distribution, consists of 1 point at the mean and $2N_{\theta}$ symmetrically distributed points around the mean. 
The forward problem evaluation at the mean is a first order approximation of the transformed mean, and the other $2N_{\theta}$ forward evaluations give a 2nd order approximation of the transformed covariance. They correspond to $\G(\theta)$ and $d \G(\theta)$ in minimizing~\cref{eq:KI2}, respectively. We observe that the covariance evaluation does not need to be very accurate, and this gives room for another level of speedup.
Reduced-order models, including but not limited to low-fidelity/low-resolution models~\cite{robinson2008surrogate,macdonald2020multi}, projection based reduced-order models (PROM)~\cite{carlberg2013gnat,qian2020lift,grimberg2020mesh} , Gaussian process based surrogate models~\cite{bilionis2012multi,cleary2020calibrate}, and neural network based surrogate models~\cite{tripathy2018deep,fan2019solving,nelsen2020random}, can be applied to speed up these $2N_{\theta}$ forward evaluations.

\subsection{Our Contributions}

Our main contribution is the development of several speedup strategies for the UKI.

\begin{itemize}
\item We take advantage of the low-rank covariance structure by reparameterization or reformulating UKI in the low-rank square root form. These strategies reduce the number of forward problem evaluations from $\bigO(N_{\theta})$ to $\bigO(N_r)$ and the memory cost from  $\bigO(N_{\theta}^2)$ to $\bigO(N_rN_{\theta})$.
\item We further speed up forward problem evaluations by using reduced-order models. The computational cost in one iteration is reduced from $2N_{\theta}+1$~($2N_{r}+1$) expensive forward evaluations to $1$ single expensive forward evaluation and $2N_{\theta}$~($2N_{r}$) cheap forward evaluations.
\item The effectiveness of the enhanced unscented Kalman inversion is demonstrated on a barotropic model inverse problem with $\bigO(10^5)$ unknown parameters and a 3D generalized circulation model, where each iteration is as efficient as that of gradient based optimization methods~(Generally, gradient based optimization methods require a forward problem evaluation and an equally expensive backward propagation in one iteration.)
\end{itemize}

The remainder of the paper is organized as follows. 
In \cref{sec:UKI}, an overview of the unscented Kalman inversion is provided.
In \cref{sec:speedup-low-rank}, speedup strategies by leveraging low-rank covariance structure are introduced. 
In \cref{sec:speedup-rom}, speedup strategies by leveraging reduced-order models are introduced. 
Numerical applications that demonstrate the efficiency of the enhanced UKI are provided in \cref{sec:app}.

\section{Unscented Kalman Inversion~(UKI)}
\label{sec:UKI}
UKI pairs the parameter-to-data map $\G$ with a stochastic dynamical system for the parameter, which is defined as following,
\begin{subequations}
  \begin{align}
  &\textrm{evolution:}    &&\theta_{n+1} = r + \alpha (\theta_{n}  - r) +  \omega_{n+1}, &&\omega_{n+1} \sim \N(0,\Sigma_{\omega}), \label{eq:dyn:evolution}\\
  &\textrm{observation:}  &&y_{n+1} = \G(\theta_{n+1}) + \nu_{n+1}, &&\nu_{n+1} \sim \N(0,\Sigma_{\nu}), \label{eq:dyn:observation}
\end{align}
\end{subequations}
where $\theta_{n+1}$ is the unknown state vector, and $y_{n+1}$ is the observation, the artificial evolution error $\omega_{n+1}$ and artificial observation error $\nu_{n+1}$ are mutually independent, zero-mean Gaussian sequences with covariances $\Sigma_{\omega}$ and $\Sigma_{\nu}$, respectively. $\alpha \in (0,1]$ is the regularization parameter,  and  $r$ is an arbitrary vector. 

Let denote  $Y_{n} := \{y_1, y_2,\cdots , y_{n}\}$, the observation set at time $n$.  
Techniques from filtering are employed to approximate the distribution $\mu_n$ of $\theta_n|Y_n$. The iterative algorithm starts from $\mu_0$ and updates 
$\mu_n$ through the prediction and
analysis steps~\cite{reich2015probabilistic,law2015data}: $\mu_n 
\mapsto \hat{\mu}_{n+1}$, and then 
$\hat{\mu}_{n+1} \mapsto \mu_{n+1}$, where $\hat{\mu}_{n+1}$ is the distribution
of $\theta_{n+1}|Y_n$. 

In the prediction step, we assume that $\mu_n \approx \N(m_n,C_n)$, then under (\ref{eq:dyn:evolution}),
$\hat{\mu}_{n+1}$ is also Gaussian with mean and covariance:
\begin{equation}
\label{eq:KF_pred_mean}
\begin{split}
    \pmean_{n+1} = \E[\theta_{n+1}|Y_n] = \alpha \mean_n + (1-\alpha)r \qquad 
    \pCov_{n+1} = \mathrm{Cov}[\theta_{n+1}|Y_n]  = \alpha^2\Cov_{n} + \Sigma_{\omega}.
\end{split}
\end{equation}
In the analysis step, we assume that the joint distribution of  $\{\theta_{n+1}, y_{n+1}\}|Y_{n}$ can be approximated by a Gaussian distribution
\begin{equation}
\label{eq:KF_joint}
     \N\Bigl(
    \begin{bmatrix}
    \pmean_{n+1}\\
    \py_{n+1}
    \end{bmatrix}, 
    \begin{bmatrix}
   \pCov_{n+1} & \pCov_{n+1}^{\theta p}\\
    {{\pCov_{n+1}}^{\theta p}}{}^{T} & \pCov_{n+1}^{pp}
    \end{bmatrix}
    \Bigr),
\end{equation}
where 
\begin{equation}
\label{eq:KF_joint2}
\begin{split}
    \py_{n+1} =     & \E[\G(\theta_{n+1})|Y_n], \\
     \pCov_{n+1}^{\theta p} =     &  \mathrm{Cov}[\theta_{n+1}, \G(\theta_{n+1})|Y_n],\\
    \pCov_{n+1}^{p p} = &  \mathrm{Cov}[\G(\theta_{n+1})|Y_n] + \Sigma_{\nu}.
\end{split}
\end{equation}
Conditioning the Gaussian in \eqref{eq:KF_joint} to find $\theta_{n+1}|\{Y_n,y_{n+1}\}=\theta_{n+1}|Y_{n+1}$ gives the following
expressions for the mean $\mean_{n+1}$ and covariance $\Cov_{n+1}$ of the
approximation to $\mu_{n+1}:$
\begin{equation}
\label{eq:KF_analysis}
    \begin{split}
        \mean_{n+1} &= \pmean_{n+1} + \pCov_{n+1}^{\theta p} (\pCov_{n+1}^{p p})^{-1} (y_{n+1} - \py_{n+1}), \\
         \Cov_{n+1} &= \pCov_{n+1} - \pCov_{n+1}^{\theta p}(\pCov_{n+1}^{p p})^{-1} {\pCov_{n+1}^{\theta p}}{}^{T}.
    \end{split}
\end{equation}
By assuming all observations $\{ y_n \}$ are identical to $y$~($Y_n = y$), \Cref{eq:KF_pred_mean,eq:KF_joint,eq:KF_joint2,eq:KF_analysis} establish a conceptual description of the Kalman inversion to solve the inverse problem~\eqref{eq:KI}.
And UKI uses the following unscented transform to evaluate \cref{eq:KF_joint2}.

\begin{definition}[Modified Unscented Transform~\cite{UKI}]
\label{def:unscented_transform}
Let denote Gaussian random variable $\theta \sim \N(\mean, \Cov) \in \R^{N_{\theta}}$, $2N_{\theta}+1$ symmetric $\sigma-$points are chosen deterministically:
\begin{align*}
    \theta^0 = \mean \qquad \theta^j = \mean + c_j [\sqrt{\Cov}]_j  \qquad \theta^{j+N_\theta}  = \mean - c_j [\sqrt{\Cov}]_j\quad (1\leq j\leq N_\theta),
\end{align*}
where $[\sqrt{\Cov}]_j$ is the $j$th column of the Cholesky factor of $\Cov$. The quadrature rule approximates the mean and covariance of the transformed variable $\G_i(\theta)$ as follows,  
\begin{equation}
\begin{split}
\label{eq:UT}
    \E[\G_i(\theta)] \approx \G_i(\theta_0) = \G_i(m)\qquad 
    \textrm{Cov}[\G_1(\theta), \G_2(\theta)]  \approx \sum_{j=1}^{2N_{\theta}} W_j^{c} (\G_1(\theta^j) - \E\G_1(\theta))(\G_2(\theta^j) - \E\G_2(\theta))^T. 
\end{split}
\end{equation}
Here these constant weights are 
\begin{align*}
    &c_1=c_2\cdots=c_{N_\theta} = \sqrt{N_\theta +\lambda} \quad W_1^{c} = W_2^{c} = \cdots = W_{2N_\theta}^{c} = \frac{1}{2(N_r+\lambda)}, \\
    &\lambda = a^2 (N_\theta + \kappa) - N_\theta \qquad \kappa = 0 \qquad a=\min\{\sqrt{\frac{4}{N_\theta + \kappa}},1\}. \end{align*}
\end{definition}

Consider the algorithm defined by
\cref{eq:KF_pred_mean,eq:KF_joint,eq:KF_joint2,eq:KF_analysis}.
By utilizing the aforementioned quadrature rule, we obtain the following UKI algorithm:
\begin{itemize}
\item Prediction step : 
    \begin{equation}
    \label{eq:UKI-pred}
        \pmean_{n+1} = \alpha \mean_{n} + (1-\alpha)r \qquad \pCov_{n+1} = \alpha^2\Cov_n + \Sigma_{\omega}.\\
    \end{equation}
    \item Generate $\sigma-$points :
    \begin{align*}
    &\pp_{n+1}^0 = \pmean_{n+1}, \\
    &\pp_{n+1}^j = \pmean_{n+1} + c_j [\sqrt{\pCov_{n+1}}]_j \quad (1\leq j\leq N_\theta),\\ 
    &\pp_{n+1}^{j+N_\theta} = \pmean_{n+1} - c_j [\sqrt{\pCov_{n+1}}]_j\quad (1\leq j\leq N_\theta).
    \end{align*}
\item Analysis step :
   \begin{subequations}
   \label{eq:UKI-analysis}
   \begin{align}
        &\ppy^j_{n+1} = \G(\pp^j_{n+1}) \qquad \py_{n+1} = \ppy^0_{n+1}, \label{eq:UKI-sample}\\
         &\pCov^{\theta p}_{n+1} = \sum_{j=1}^{2N_\theta}W_j^{c}
        (\pp^j_{n+1} - \pmean_{n+1} )(\ppy^j_{n+1} - \py_{n+1})^T, \\
        &\pCov^{pp}_{n+1} = \sum_{j=1}^{2N_\theta}W_j^{c}
        (\ppy^j_{n+1} - \py_{n+1} )(\ppy^j_{n+1} - \py_{n+1})^T + \Sigma_{\nu},\\
        &\mean_{n+1} = \pmean_{n+1} + \pCov^{\theta p}_{n+1}(\pCov^{pp}_{n+1})^{-1}(y - \py_{n+1}),\label{eq:uki-mean-prop}\\
        &\Cov_{n+1} = \pCov_{n+1} - \pCov^{\theta p}_{n+1}(\pCov^{pp}_{n+1})^{-1}{\pCov^{\theta p}_{n+1}}{}^{T}.\label{eq:uki-cov-prop}
      \end{align}
       \end{subequations}
\end{itemize}
Following~\cite{UKI}, the hyperparameters of the stochastic dynamical system are 
\begin{equation}
\label{eq:hyperparameters}
   r = r_0, \quad \Sigma_{\omega} = (2 - \alpha^2)\Lambda, \quad \Sigma_{\nu} = 2\Sigma_{\eta}, \quad \textrm{ and } \quad \alpha \in (0, 1.0]. 
\end{equation}
where $r_0$ is the prior mean and $\Lambda$ is a positive definite matrix. And the initial distribution $\theta_0\sim \N(\mean_0, \Cov_0)$ satisfies $\mean_0 = r_0$ and $\Cov_0 = \Lambda$.

UKI~\cite{UKI} as a general derivative-free approach to solving the inverse problem, is particularly suitable for inverse problems where the forward model  is  given  as  a  black  box  and  may  not  be  differentiable.
However, the computational cost of $2N_\theta+1$ forward model evaluations~\eqref{eq:UKI-sample} in the analysis step could be intractable, especially when the unknown parameter dimension is high or the forward problem evaluation is expensive; the memory cost of storing and propagating the covariance~\eqref{eq:uki-cov-prop} could be unaffordable, especially when the unknown parameter dimension is high. In the present work, we focus on improving UKI by leveraging the low rank covariance structure and the reduced-order model techniques.

\section{Speed Up with Low-Rank Covariance Structure}
\label{sec:speedup-low-rank}
In general, UKI as a particle method, requires evaluations on $2N_{\theta}+1$ $\sigma$-points and the propagation of the $N_\theta \times N_\theta$ covariance matrix, which is impractical for general geophysical applications~(i.e., initial condition recovery problems) with dimension of roughly $\bigO(10^4) - \bigO(10^6)$~\cite{van2003variance,oliver2011recent}. However, the ensemble Kalman filter~\cite{evensen1994sequential,anderson2001ensemble,bishop2001adaptive} is successfully used in high-dimensional data-assimilation applications with ensemble size $\bigO(10^2)$, what is the catch? 
Since many observed geophysical phenomena have large-scale structure or low-dimensional attractors, therefore, 
the model error covariance can be well-approximated by the low-rank sample covariance.
Moreover, these large-scale structure can be well represented by dominant frequency modes~(i.e., dominant Fourier or Karhunen-Lo\`{e}ve~(KL) modes).

In this section, two approaches to leveraging the low-rank large-scale structure in the UKI framework are introduced: 
\begin{itemize}
    \item Reparameterization, namely project or reparameterize the model parameters onto some vastly reduced subspace;
    \item Truncated unscented Kalman inversion, namely formulate the unscented Kalman filter as a low-rank square-root filter.
\end{itemize}
For both approaches, the number of $\sigma$-points is reduced from $2N_{\theta}+1$ to $2N_r+1$, where $N_r$ is the low rank number. And the storage for the covariance matrix is reduced from $\bigO(N_{\theta}^2)$ to $\bigO(N_{\theta} N_r)$. Moreover, a discussion about the scenario, where low rank covariance structure is unknown or does not exist, is presented in Subsection~\ref{ssec:no-low-rank}.

\subsection{Reparameterization}
Let $r_0$ denote the prior mean. The discrepancy $\theta - r_0$ is assumed to be well-approximated in the linear space spanned by the basis $\{u_1, u_2,\cdots, u_{N_r}\}$. Hence, the unknown parameters can be reparameterized as following, 
\begin{equation*}
    \theta = r_0 + \sum_{i=1}^{N_r}\tau_{(i)} u_i.
\end{equation*}
UKI~\cite{UKI} is then applied to solve for the vector $\tau = [\tau_{(1)}, \tau_{(2)}, \cdots, \tau_{(N_r)}]^T$, which has prior mean $0$.

Let $\mean_n(\tau)$ and $\Cov_n(\tau)$ denote the mean and covariance estimation of $\tau$, the mean and covariance estimation of $\theta$ can be recovered by
\begin{equation}
\label{eq:reparam_mean_cov}
  \mean_n(\theta) = r_0 + U \mean_n(\tau) \qquad  \Cov_n(\theta) = U \Cov_n(\tau) U^T,
\end{equation}
here $\displaystyle U = \Big(u_1\,\, u_2\,\cdots\, u_{N_r}\Big)$.

\begin{remark}
\label{rem:reparam}
\Cref{eq:reparam_mean_cov} also reveals the low-rank covariance structure of $\Cov_n(\theta)$.
The square root\footnote{
In the present work, the square root of the matrix $C$ is defined to be a matrix $Z$ such that 
\begin{equation*}
    C = Z Z^T,
\end{equation*}
which is inconsistent with its most common use in the mathematical literature.
And if $Z_1$ and $Z_2$ are two $n\times m$ square root matrices of $C$, then there exists an orthogonal matrix $Q$ such that
$Z_2 = Z_1 Q$~\cite{livings2008unbiased}.
}
$Z_n$ of the covariance matrix $\Cov_n(\theta)$ can be written as 
\begin{equation*}
    Z_n = U \sqrt{\Cov_n(\tau)} Q,
\end{equation*}
here $Q$ is an orthogonal matrix. Hence, all the square root matrices $\{Z_n\}$ share the column vector space---the space spanned by 
column vectors of $U$, which is prescribed in the parameterization.
\end{remark}

\subsection{Truncated Unscented Kalman Inversion~(TUKI)}
\label{subsec:UKI}
The truncated unscented transform, which is a low-rank approximation of the modified unscented transform~(See \cref{def:unscented_transform}), is central to the TUKI:
\begin{definition}[Truncated Unscented Transform]
\label{def:truncated_unscented_transform}
Let denote Gaussian random variable $\theta \sim \N(\mean, \Cov) \in \R^{N_{\theta}}$,
and assume the $N_r$ dimensional truncated singular value decomposition~($N_r$-TSVD) of $C$ is
\begin{align*}
    C \approx \sum_{i=1}^{N_{r}}d_i u_i u_i^T,
\end{align*}
where $\{u_i\}$ are left singular vectors, and $d_1 \geq d_2  \cdots  d_{N_{r}}\geq 0$ are dominant singular values.
$2N_{r}+1$ symmetric $\sigma-$points are chosen deterministically:
\begin{align*}
    \theta^0 = \mean \qquad \theta^j = \mean + c_j \sqrt{d_j} u_j  \qquad \theta^{j+N_r}  = \mean - c_j \sqrt{d_j} u_j \quad (1\leq j\leq N_r).
\end{align*}
The quadrature rule approximates the mean and covariance of the transformed variable $\G_i(\theta)$ as follows,  
\begin{align*}
    \E[\G_i(\theta)] \approx  \G_i(\theta^0)\qquad 
    \textrm{Cov}[\G_1(\theta),\G_2(\theta)]  \approx \sum_{j=1}^{2N_r} W_j^{c} (\G_1(\theta^j) - \E\G_1(\theta))(\G_2(\theta^j) - \E\G_2(\theta))^T. 
\end{align*}
Here these constant weights are 
\begin{align*}
    &c_1=c_2\cdots=c_{N_r} = \sqrt{N_r +\lambda} \quad W_1^{c} = W_2^{c} = \cdots = W_{2N_r}^{c} = \frac{1}{2(N_r+\lambda)},\\
    &\lambda = a^2 (N_r + \kappa) - N_r \qquad \kappa = 0 \qquad a=\min\{\sqrt{\frac{4}{N_r + \kappa}},1\}.
\end{align*}
\end{definition}

We choose the hyperparameter $\Lambda$ in~\cref{eq:hyperparameters} to be a low-rank matrix with square root $Z_0 \in R^{N_\theta \times N_r}$, which satisfies
\begin{equation}
\label{eq:TUKI-start}
    \Lambda = Z_0 Z_0^T \quad \textrm{and} \quad \textrm{rank}(Z_0) = N_r.
\end{equation}
Let $Z_{\omega} = \sqrt{2 -\alpha^2}Z_0$ denote the square root matrix of the artificial evolution error covariance $\Sigma_{\omega}$ in~\cref{eq:hyperparameters}.
The iteration procedure of the UKI can be formulated in the square root form $\{Z_n\}$, where $\Cov_n = \CovZ_n \CovZ_n^T$.
This leads to the following truncated unscented Kalman inversion~(TUKI):
\begin{itemize}
\item Prediction step : 
    \begin{equation}
    \label{eq:TUKI-pred}
        \pmean_{n+1} = \alpha \mean_{n} + (1-\alpha)r \qquad N_r\textrm{-TSVD : }\Big(\alpha\CovZ_n \quad \CovZ_{\omega}\Big) = \widehat{U}_n \sqrt{\widehat{D}_n} \widehat{V}_n^T,
    \end{equation}
    where $\displaystyle \widehat{U}_n = \Big(u_1\,\, u_2\,\cdots\, u_{N_r}\Big)$ and $\widehat{D}_n = \textrm{diag}\{d_1, d_2 ,\cdots d_{N_r}\}$. And hence
    $$\pCov_{n+1} = \alpha^2\Cov_n + \Sigma_{\omega} =
    \Big(\alpha\CovZ_n \quad \CovZ_{\omega}\Big)\Big(\alpha\CovZ_n \quad \CovZ_{\omega}\Big)^T = \widehat{U}_n \widehat{D}_n \widehat{U}_n^T.$$
    \item Generate truncated $\sigma-$points :
    \begin{align*}
    &\pp_{n+1}^0 = \pmean_{n+1} \qquad \pp_{n+1}^j = \pmean_{n+1} + c_j \sqrt{d_j} u_j \quad \pp_{n+1}^{j+N_r} = \pmean_{n+1} - c_j \sqrt{d_j} u_j\quad (1\leq j\leq N_r).
    \end{align*}
\item Analysis step :
   \begin{subequations}
   \label{eq:TUKI-analysis}
   \begin{align}
        &\ppy^j_{n+1} = \G(\pp^j_{n+1}) \qquad \py_{n+1} = \ppy^0_{n+1},\label{eq:TUKI-sample}\\
        &\pCovZ_{n+1} = \Big(\sqrt{W_1^c}(\pp^1_{n+1} - \pmean_{n+1} )\quad \sqrt{W_2^c}(\pp^2_{n+1} - \pmean_{n+1} )\quad\cdots\quad \sqrt{W_{2N_r}^c}(\pp^{2N_r}_{n+1} - \pmean_{n+1} ) \Big),\\
        &\pGCovZ_{n+1} = \Big(\sqrt{W_1^c}(\ppy^1_{n+1} - \ppy_{n+1} )\quad \sqrt{W_2^c}(\ppy^2_{n+1} - \ppy_{n+1} )\quad\cdots\quad \sqrt{W_{2N_r}^c}(\ppy^{2N_r}_{n+1} - \ppy_{n+1} ) \Big), \\
        &\textrm{SVD : }\pGCovZ_{n+1}^T \Sigma_{\nu}^{-1}\pGCovZ_{n+1} = \widehat{P}_{n+1} \widehat{\Gamma}_{n+1} \widehat{P}_{n+1}^T,\label{eq:TUKI-update-svd}\\
         &\mean_{n+1} = \pmean_{n+1} + \pCovZ_{n+1}\widehat{P}_{n+1}(\widehat{\Gamma}_{n+1}+\I)^{-1}\widehat{P}_{n+1}^T
        \pGCovZ_{n+1}^T\Sigma_{\nu}^{-1}(y - \py_{n+1}),\label{eq:TUKI-update-mean}\\
        &\CovZ_{n+1} = \pCovZ_{n+1}\widehat{P}_{n+1}(\widehat{\Gamma}_{n+1}+\I)^{-1/2}.\label{eq:TUKI-update-cov}
    \end{align}
    \end{subequations}
\end{itemize}

Since the covariance matrices in~\cref{eq:UKI-analysis} are    %
\begin{equation*}
 \pCov^{\theta p}_{n+1} = \pCovZ_{n+1} \pGCovZ_{n+1}^T, \quad 
\pCov^{pp}_{n+1} = \pGCovZ_{n+1} \pGCovZ_{n+1}^T + \Sigma_{\nu}, \quad \textrm{ and }
\pCov_{n+1} = \pCovZ_{n+1} \pCovZ_{n+1}^T. 
\end{equation*}
Applying Sherman–Morrison–Woodbury formula
\begin{equation*}
    (\pGCovZ_{n+1} \pGCovZ_{n+1}^T + \Sigma_{\nu})^{-1} = \Sigma_{\nu}^{-1} - \Sigma_{\nu}^{-1}\pGCovZ_{n+1}(\I + \pGCovZ_{n+1}^T \Sigma_{\nu}^{-1}\pGCovZ_{n+1})^{-1} \pGCovZ^T_{n+1}\Sigma_{\nu}^{-1}
\end{equation*}
to these update equations~\eqref{eq:uki-mean-prop} and \eqref{eq:uki-cov-prop} leads to the update equations \eqref{eq:TUKI-update-mean} and \eqref{eq:TUKI-update-cov}.

\begin{remark}
\label{rem:SVD}
The SVD in~\cref{eq:TUKI-update-svd} is applied on a $2N_r \times 2N_r$ matrix and the truncated SVD in~\cref{eq:TUKI-pred} is applied on a $N_\theta \times 2N_r$ matrix, and hence, they are affordable.
\end{remark}

\begin{lemma}
\label{lemma:TUKI}
All square root matrices $\{\CovZ_n\}$ and $\{\pCovZ_n\}$ are rank $N_r$ and share the same column vector space spanned by the column vectors of $Z_0$~(\ref{eq:TUKI-start}). Therefore, the $N_r$-TSVD in the prediction step~\eqref{eq:TUKI-pred} is exact. 
\end{lemma}
\begin{proof}
The proof is in \ref{sec:app:proof}.
\end{proof}

\begin{remark}
\label{rem:invariant}
$\{m_{n}\}$ lie in the linear space spanned by $r$~(or $r_0$) and columns of $Z_0$, since 
$$\mean_{n+1} = \alpha\mean_n  + (1-\alpha)r  + \pCovZ_{n+1}\widehat{P}_{n+1}(\widehat{\Gamma}_{n+1}+\I)^{-1}\widehat{P}_{n+1}^T
\pGCovZ_{n+1}^T\Sigma_{\nu}^{-1}(y - \py_{n+1}).
$$
\end{remark}
\begin{remark}
Although the square root matrix $\CovZ_{n+1}$ computed in~\cref{eq:TUKI-update-cov} is not unique,  namely it can be represented as 
\begin{equation*}
\CovZ_{n+1} = \pCovZ_{n+1}\widehat{P}_{n+1}(\widehat{\Gamma}_{n+1}+\I)^{-1/2} Q,
\end{equation*}
where $Q$ can be any orthogonal matrix, the algorithm is uniquely defined. Since the SVD of the square root of the predicted covariance matrix~\eqref{eq:TUKI-pred} can be written as  
\begin{equation*}
    \Big(\alpha\CovZ_{n+1} Q \quad \CovZ_{\omega}\Big) = \widehat{U}_{n+1} \sqrt{\widehat{D}_{n+1}} \Big( \widehat{V}_{n+1,1}Q \quad \widehat{V}_{n+1,2} \Big),
\end{equation*}
where $ \displaystyle \Big(\alpha\CovZ_{n+1} \quad \CovZ_{\omega}\Big) = \widehat{U}_{n+1} \sqrt{\widehat{D}_{n+1}} \Big( \widehat{V}_{n+1,1} \quad \widehat{V}_{n+1,2} \Big) $, 
the orthogonal matrix $Q$ does not affect the left singular vectors $\widehat{U}_{n+1}$ or singular values in $\widehat{D}_{n+1}$.
\end{remark}
\subsection{Without Low Rank Covariance Structure Information}
\label{ssec:no-low-rank}
In this section, we focus on problems with high dimensional unknown parameters, but without low-rank covariance structure information~(i.e., the prior covariance is $\I$).
For ensemble Kalman inversion and its variants~(See \ref{sec:app:ensemble_method}), the initial ensemble can still be generated from the prior distribution. 
For TUKI, a low rank matrix $\Lambda = Z_0 Z_0^T$ is required to reduce computational and memory costs.

However, the ensemble Kalman inversion has the invariant subspace property~\cite{iglesias2013ensemble,schillings2017analysis}, namely, the ensemble members are in the linear space spanned by the initial ensembles. 
The UKI with reparameterization and TUKI have  similar properties~(See~\cref{rem:reparam,rem:invariant}). 
On one hand, this property brings regularization effect. On the other hand, this property also brings strong limitation for the inverse problems without low-rank covariance information, when the truncated low rank number $N_r$ or the ensemble size $J$  are much smaller than the unknown parameter dimension~$N_{\theta}$. The limitation is illustrated in the following theorem:
\begin{theorem}
\label{theorem:negative}
Let $\{\theta^i\}_{i=1}^J$ be random samples from the $N_\theta$-dimensional Gaussian distribution $\N(0, \I)$, for any given vector $\theta_{ref} \in \R^{N_{\theta}}$, we have 
\begin{equation}\label{e:lob}
    \E \bigl[\textrm{dist}(\theta_{ref}, \textrm{span}\{\theta^1, \theta^2, \cdots, \theta^J\})^2 \bigr] =(1-\frac{J}{N_\theta})\|\theta_{ref}\|^2.
\end{equation}
\end{theorem}
\begin{proof}
The proof is in \ref{sec:app:proof}.
\end{proof}

\cref{theorem:negative} illustrates that when $J \textrm { and } N_r$ are much smaller than $N_{\theta}$, even the optimal approximation in the randomly generated invariant subspace is a poor approximation of the solution with high probability. 

\section{Speed Up with Reduced-Order Models}
\label{sec:speedup-rom}
For large scientific or engineering problems, even with a small parameter number $N_\theta$~(or rank number $N_r$), the computational cost associated with these $2N_{\theta}+1$~(or $2N_{r}+1$) forward model evaluations can be intractable. Reduced-order models can be used to speed up these $2N_{\theta}+1$~(or $2N_{r}+1$) forward model evaluations.
Specifically the forward evaluation $\G(\hat{\theta}_{n+1}^0)$ is evaluated by the high-fidelity model, and other forward evaluations
\begin{equation*}
    \G(\hat{\theta}_{n+1}^j)\qquad  \textrm{with} \qquad j \geq 1
\end{equation*}
in \cref{eq:UKI-sample,eq:TUKI-sample} are evaluated with reduced-order models.
Reduced-order models include but not limited to low-fidelity/low-resolution models~\cite{robinson2008surrogate,macdonald2020multi}, projection based reduced-order models (PROM)~\cite{qian2020lift,carlberg2013gnat,grimberg2020mesh}, Gaussian process based surrogate models~\cite{bilionis2012multi,cleary2020calibrate}, and neural network based surrogate models~\cite{tripathy2018deep,fan2019solving,nelsen2020random}.
For low-fidelity/low-resolution models, they can be applied straightforwardly.
For PROMs, a training procedure is required. The data generated at $\pp_{n+1}^0$ can form a good projection basis. For Gaussian process and neural network based surrogate models, data generating and training procedures are required.

\begin{remark}
Following Proposition 2 in~\cite{UKI}, the UKI~($\alpha=1$) update equations~(\ref{eq:UKI-analysis}) can be written as 
\begin{subequations}
\begin{align}
&\mean_{n+1} = \mean_{n} + \Cov_{n+1}\F_u d\G_{n+1}^T\big(\Sigma_{\nu} + \widetilde{\Sigma}_{\nu,n+1}\big)^{-1}\big(y - \G(\mean_n)\big)\label{eq:uki-mean},\\
&\Cov_{n+1} = \Big(\pCov_{n+1}^{-1} + \F_u d\G_{n+1}^T \big(\Sigma_{\nu} + \widetilde{\Sigma}_{\nu,n+1}\big)^{-1} \F_u d\G_{n+1}\Big)^{-1}.
\end{align}
\end{subequations}
where $\F_u d\G_{n+1} = {\pCov_{n+1}^{\theta p}}{}^T {\pCov_{n+1}}^{-1}$ is the averaged derivative and $\|\widetilde{\Sigma}_{\nu,n+1}\| = \bigO(\|\pCov_{n+1}^2\|)$.
%
To evaluate $\G(\mean_n)$ and $\F_u d\G_{n+1}$, the UKI requires $2N_{\theta}+1$ forward model evaluations. Specifically, in the modified unscented transform~(See \cref{def:unscented_transform}), we need $1$ forward model evaluation 
$$ \py_{n+1}^0 = \G(\pp_{n+1}^0) = \G(\mean_{n}) $$
to estimate the mean, which corresponds to $y - \G(\mean_n)$ in \cref{eq:uki-mean}; and $2N_\theta$ forward model evaluations 
$$ \py_{n+1}^j = \G(\pp_{n+1}^j) \qquad j \geq 1$$
to estimate the covariance, which corresponds to  $\F_u d\G_{n+1}$. 

Since we find the covariance and $\F_u d\G_{n+1}$ do not need to be very accurate. Hence, these $2N_\theta$  forward model evaluations can be accelerated by using reduced-order models.
\end{remark}

\section{Applications}
\label{sec:app}
In this section, we present numerical results for using the aforementioned strategies, including UKI with  reparameterization~(referred as UKI) and TUKI, to speed up UKI. 
For comparison, we also apply ensemble Kalman inversion~(EKI), ensemble adjustment Kalman inversion~(EAKI), and ensemble transform Kalman inversion~(ETKI), the detailed implementation of these ensemble-based Kalman inversions is listed in \ref{sec:app:ensemble_method}.

\begin{itemize}
    \item Linear model problems: these problems serve as proof-of-concept examples, which compare the behaviors of different Kalman inversions with/without low-rank covariance structure information.
    \item Barotropic model problem: this is a model data assimilation problem in meteorology, where the initial condition is recovered from random observations in the north hemisphere. This problem demonstrates the effectiveness of using low rank covariance structure information to speed up UKI.
    \item Idealized general circulation model problem: this is a 3D Navier-Stokes problem with hydrostatic assumption, which describes physical processes in the atmosphere. This problem demonstrates the effectiveness of using reduced-order models~(specifically low-resolution models) to speed up UKI. 
\end{itemize}

The code is accessible online:
\begin{center}
  \url{https://github.com/Zhengyu-Huang/InverseProblems.jl}
\end{center}

\subsection{Linear Model Problems}
In this section, 2 linear model problems are considered. The first one is from one dimensional elliptic equation and, therefore the low-rank covariance structure exists. The second one is an artificial linear model problem, with parameters following Bernoulli distribution and, therefore the low-rank covariance structure information does not exist.
For comparison, UKI with reparameterization, TUKI, EKI, EAKI, and ETKI are all applied to these linear problems. Since these problems are well-posed and no observation error is added, the regularization parameter $\alpha$ is set to be 1.

\subsubsection{With Low Rank Covariance Structure}
\label{sec:elliptic}
Consider the one dimensional elliptic equation
\begin{equation*}
\begin{split}
&-\frac{d^2 \theta}{dx^2} + \theta  = f(x), \\
&\theta(0) = \theta(1) = 0.
\end{split}
\end{equation*}
Here  Dirichlet boundary conditions are applied on both ends, and $f$ defines the source:
\begin{align*}
    f(x) = \begin{cases}
               1 & 0 \leq x\leq \frac{1}{2}\\
               2 & \frac{1}{2} < x \leq 1\\
            \end{cases}. 
\end{align*}

The elliptic equation is semi-discretized by finite difference method on a uniform grid with $N_{\theta}+2$ points with $N_{\theta} = 1000$. 
The solution is $\theta_{ref} \in \R^{N_{\theta}}$ on these interior points.
We are interested in the inverse problem of recovering $\theta_{ref}$ from observations of $f$ on these interior points ($N_y = N_\theta$):
\begin{equation*}
    y =  G\theta + \eta,
\end{equation*}
here $\displaystyle G = (-\frac{d^2}{dx^2} + 1)$ is the discretized operator, the observation error is $\eta \sim \N(0, \I) \in \R^{N_y}$. 
Since the large structure of $\theta$ lies in the Fourier sine space, we can approximate it as 
\begin{equation}
\label{eq:elliptic-reparam}
    \theta = \tau_{(1)}\sin(x\pi)+\,\cdots\, + \tau_{(N_r)} \sin(N_rx\pi) \quad \textrm{  with  } \quad N_r = 5.
\end{equation}
For UKI, the reparameterization is \cref{eq:elliptic-reparam} and the UKI is initialized with $\tau_0 \sim \N(0,10^2\I)$. 
For TUKI, the low-rank covariance structure is embedded in the square root of $\Lambda$:
\begin{equation*}
    Z_0 = 10\times\Big(\sin(x\pi)\quad \sin(2x\pi)\quad\cdots \quad\sin(N_rx\pi)\Big).
\end{equation*}
And the TUKI is initialized with $\theta_0 \sim \N(0, Z_0 Z_0^T)$. 
For EKI, EAKI, and ETKI, the initial ensembles are generated from the distribution $\N(0, Z_0 Z_0^T)$, which span the column space of $Z_0$. The ensemble size is $J= 2N_r+1$, which matches the number of $\sigma$-points of UKI and TUKI.

\begin{figure}[ht]
\centering
\includegraphics[width=0.5\textwidth]{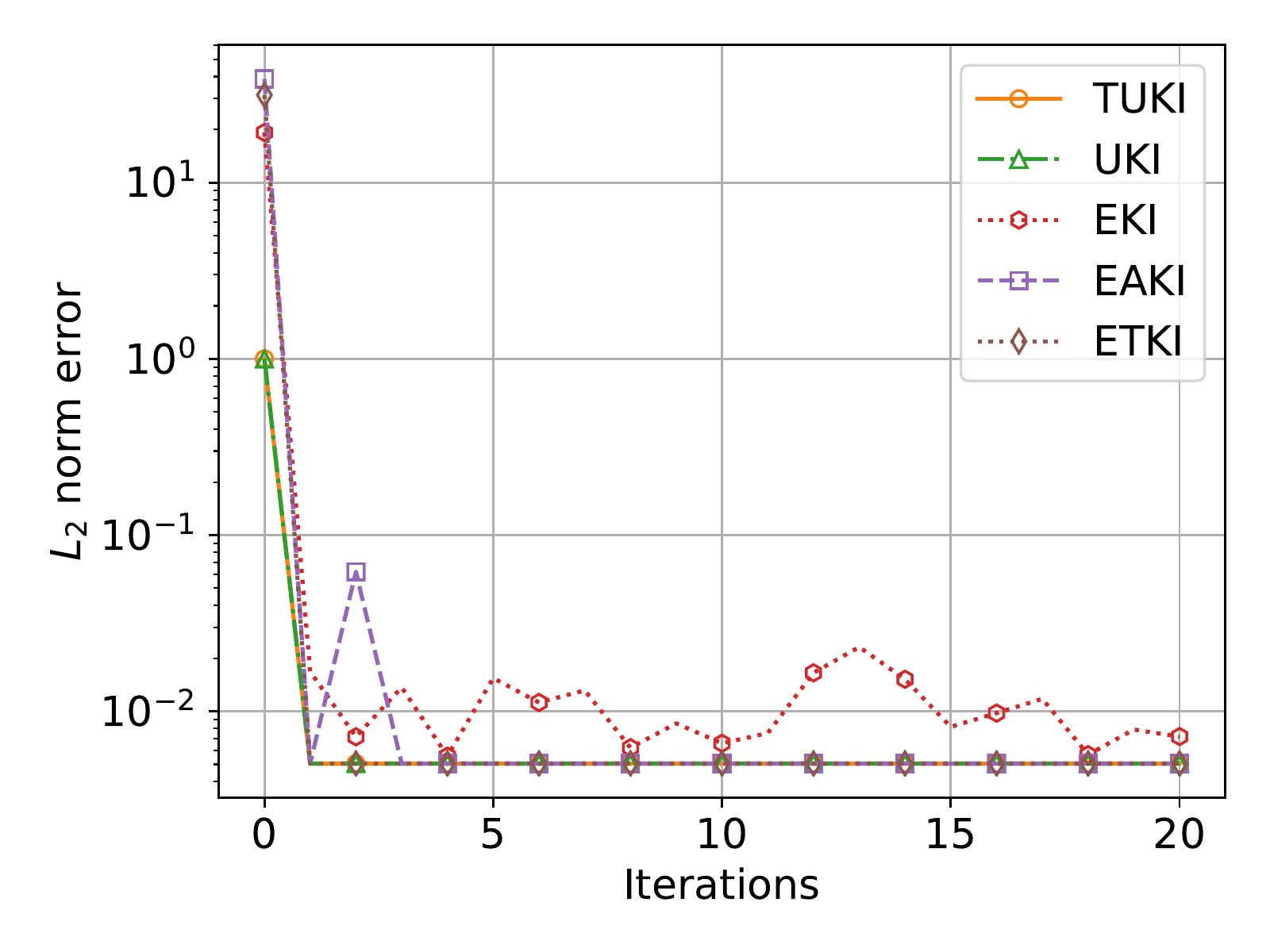}
\caption{$L_2$ error $\frac{\lVert\mean_n - \theta_{ref}\rVert_2}{\Vert\theta_{ref}\rVert_2}$ of the elliptic linear model problem.}
\label{fig:Low-rank-linear}
\end{figure}

The convergence of $\{\mean_n\}$ in terms of the relative $L_2$ norm errors is depicted in \cref{fig:Low-rank-linear}. All the Kalman inversions behave similar and converge efficiently, although the EKI suffers from random noise due to the small ensemble size. This problem illustrates how UKI uses $2N_r+1 = 11$ $\sigma$-points to solve a $N_{\theta}=1000$ dimensional inverse problem.

\subsubsection{Without Low Rank Covariance Structure}

Consider the following artificial inverse problem 
\begin{equation*}
\begin{split}
y = G\theta + \eta, \\
\end{split}
\end{equation*}
with $G = \I\in\R^{N_{y} \times N_{\theta}}$ and $N_{y} = 1000\,,\,N_{\theta} = 1000$. 
The entries of the underlying truth $\theta_{ref}$ are independent, identical Bernoulli distribution $B(1, 0.5)$. The observation is $y = G\theta_{ref}$, with observation error $\eta \sim \N(0,\I)$. 
However, there is no low-rank covariance structure information. Hence, the initial distribution is set to be $\N(0,  \I)$. 

Although there is no low-rank covariance structure information, the TUKI is applied, with truncated rank $N_r = 5$.
EKI, EAKI, and ETKI are also applied with the ensemble size of $J= 2N_r+1$.

\begin{figure}[ht]
\centering
\includegraphics[width=0.5\textwidth]{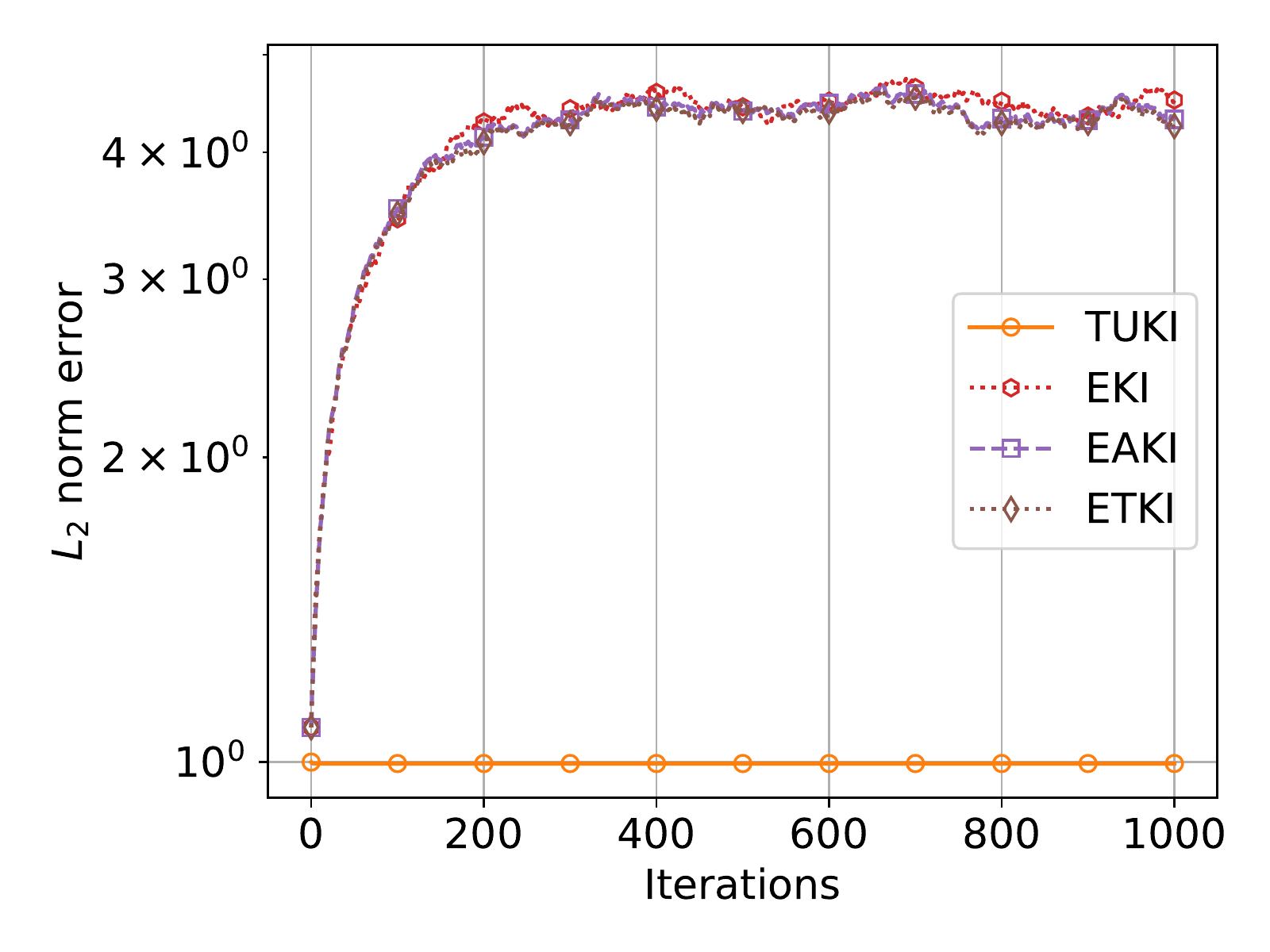}
\caption{$L_2$ error $\frac{\lVert\mean_n - \theta_{ref}\rVert_2}{\Vert\theta_{ref}\rVert_2}$ of the Bernoulli linear model problem.}
\label{fig:No-low-rank-linear}
\end{figure}

The convergence of $\{\mean_n\}$ in terms of the relative $L_2$ norm errors is depicted in \cref{fig:No-low-rank-linear}. All EKI, EAKI, and ETKI diverge, and TUKI fails to converge to $\theta_{ref}$. This conforms to the discussion in Subsection~\ref{ssec:no-low-rank}.
Gradient-based optimization algorithms with automatic differentiation~\cite{sambridge2007automatic,chen2018taylor,huang2020learning,xu2021learning} could be much more effective and efficient for this scenario, when the forward problem is differentiable and not chaotic.

\subsection{Barotropic Model Problem}
The barotropic vorticity equation describes the evolution of a non-divergent, incompressible flow on the surface of the earth:
\begin{equation}
    \begin{split}
        &\frac{\partial \omega}{\partial t} = - v \cdot \nabla (\omega + f), \\
        &\nabla^2\psi = \omega \qquad v = k \times \nabla\psi,
    \end{split}
\end{equation}
where $\omega$ and $\psi$ are (absolute) vorticity and streamfunction, respectively. $v$ is the non-divergent flow velocity, $k$ is the unit vector in the radial direction, and $f = 2\Omega \sin(\theta)$ is the Coriolis force, depending on the latitude $\theta$. The angular velocity is $\Omega = 7.292\times10^5 s^{-1}$ and the Earth radius is $R = 6.3712\times10^6~m$.

\begin{figure}[ht]
    \centering
    \begin{subfigure}{0.49\textwidth}
    \centering
        \includegraphics[width=0.9\linewidth]{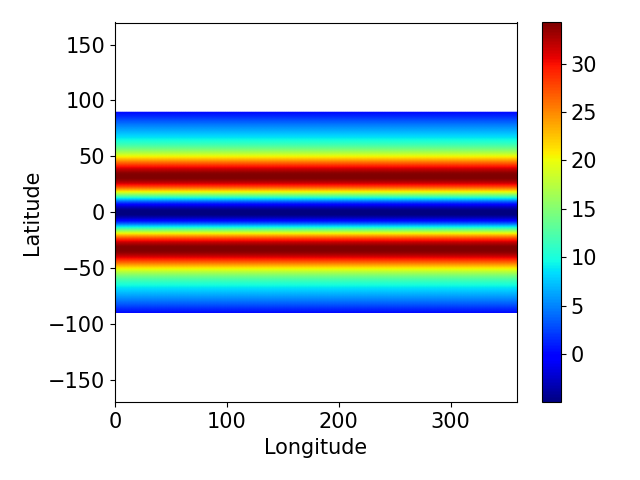}
        \caption{Zonally symmetric flow~($u_b$)}
    \end{subfigure}%
    ~ 
    \begin{subfigure}{0.49\textwidth}
    \centering
        \includegraphics[width=0.9\linewidth]{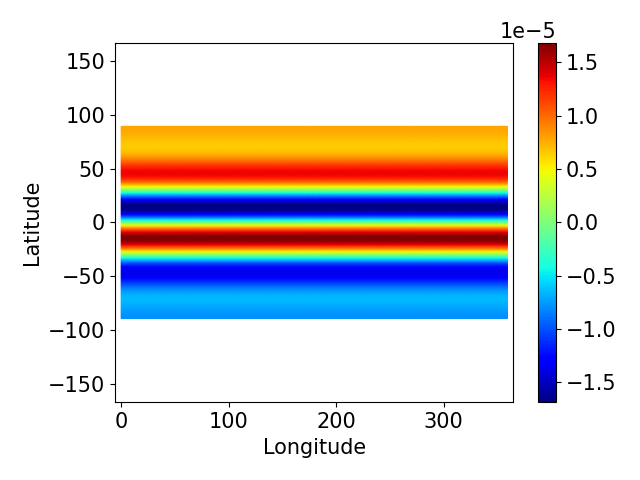}
        \caption{Zonally symmetric flow~($\omega_b$)}
    \end{subfigure}%
    
    \begin{subfigure}{0.49\textwidth}
    \centering
        \includegraphics[width=0.9\linewidth]{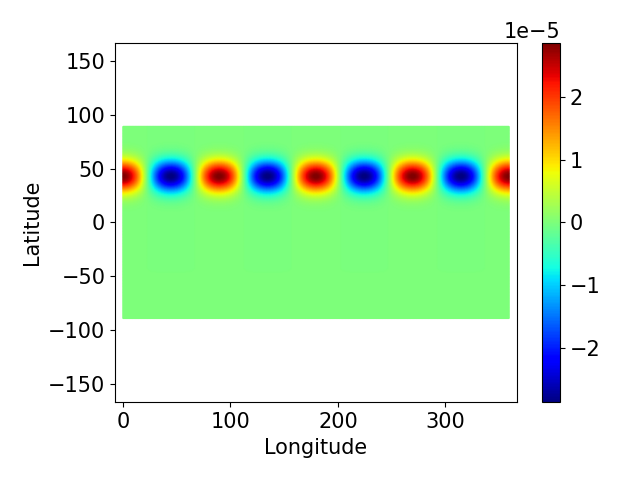}
        \caption{Sinusoidal vorticity disturbance~($\omega'$) }
    \end{subfigure}%
    ~ 
    \begin{subfigure}{0.49\textwidth}
    \centering
        \includegraphics[width=0.9\linewidth]{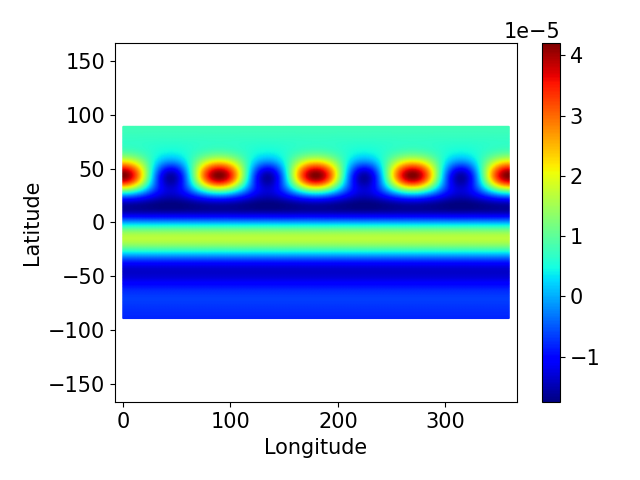}
        \caption{Initial vorticity field~($\omega_0$)}
    \end{subfigure}%
    \caption{Initial condition~(d) of the barotropic model is the superposition of a zonally symmetric flow~(a or b) and a sinusoidal disturbance in the vorticity field~(c).}
    \label{fig:vor_0}
\end{figure}

The initial condition~(See \cref{fig:vor_0}) is the superposition of a zonally symmetric flow
\begin{equation}
\label{eq:u_b}
    u_b = 25 \cos(\theta) - 30 \cos^3(\theta) + 300 \sin^2(\theta)\cos^6(\theta),
\end{equation} 
here $u_b$ is the zonal wind velocity, and a sinusoidal disturbance in the vorticity field
\begin{equation*}
    \label{eq:vor_disturbance}
    \omega' = \frac{A}{2} \cos(\theta)e^{- ((\theta - \theta_0)/\theta_W )^2} \cos(m \lambda), 
\end{equation*}
where $\lambda$ is the longitude, $m = 4, \theta_0 = 45^{\circ} N, \theta_W = 15^{\circ}$, and $A = 8.0\times10^{-5}$.
This resembles that found in the upper troposphere in Northern winter~\cite{held1985pseudomomentum,held1987linear}.

The forward problem is solved by the spectral transform method with T85 spectral resolution~(triangular truncation at wavenumber 85, with $256 \times 512$ points on the latitude-longitude transform grid) and the semi-implicit time integrator~\cite{held1994proposal}. 

\begin{figure}[ht]
\centering
\includegraphics[width=0.48\textwidth]{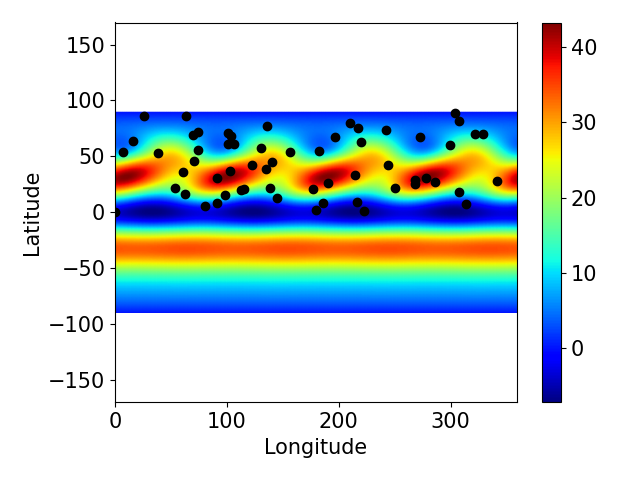}
\includegraphics[width=0.48\textwidth]{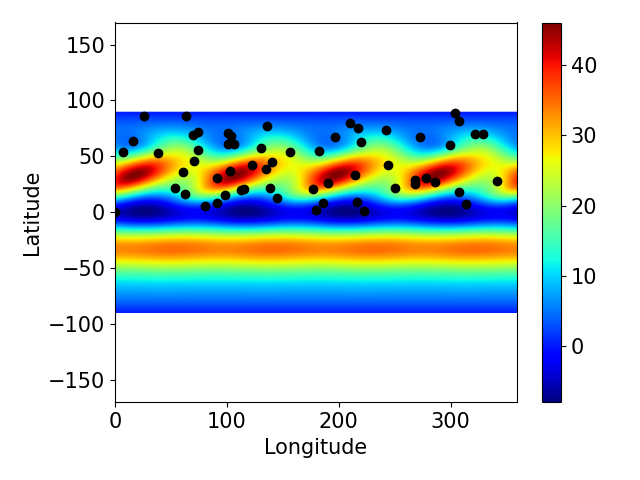}
\caption{The zonal velocity field of the barotropic model and the $50$ random pointwise measurements~(black dots) in the north hemisphere at two observation times~($T=12h$ and $T=24h$).}
\label{fig:barotropic-obs}
\end{figure}

For the inverse problem, we recover the initial condition, specifically the initial vorticity field~($\omega_0$) of the barotropic vorticity equation, given noisy pointwise observations of the zonal wind velocity field. The pointwise observations are collected at 50 random points in the north hemisphere at $T=12h$ and $T=24h$, therefore, we have $N_y=100$ data~(See \cref{fig:barotropic-obs}). 
And $5\%$ Gaussian random noises are added to the observation, as follows,
\begin{equation*}
    y_{obs} = y_{ref} + 5\% y_{ref} \odot \N(0, \I),
\end{equation*}
here $\odot$ denotes element-wise multiplication.
The initial zonal flow $\omega_b$ represents the basic atmospheric circulation, and we are more interested in the initial perturbation field $\omega'$. Hence, we assume the prior mean is $r_0 = \omega_b$.
Since the large structure of $\omega_0$ lies in the spherical harmonics space
\begin{equation*}
    Y_{m,n}(\lambda,\theta) = P_{m,n}( \theta) e^{im\lambda}\qquad  -n \leq m \leq n,\, n\geq 0.
\end{equation*}
The vorticity field on the sphere can be approximated with $N = 7$ truncated wavenumber with triangular truncation, as following,
\begin{equation}
\label{eq:sph-harmonics}
    \omega_0(\lambda, \theta) = \omega_b + \sum_{n=1}^{N} \tau_{(0,n)} P_{0,n}(\theta) + \sum_{m = 1}^{N} \sum_{n=m}^{N} \tau^s_{(m,n)} 2P_{m,n}(\theta) \cos(m\lambda) +  \tau^c_{(m,n)}2P_{m,n}(\theta) \sin(m\lambda). \\
\end{equation}
For UKI, the reparameterization is \cref{eq:sph-harmonics}, and the UKI is initialized with $\displaystyle \tau_0 \sim \N(0, \frac{\I}{R^2})$. 
For TUKI, we assume the square root of $\Lambda$ is
\begin{equation*}
    Z_0 = \frac{1}{R}\Big(P_{0,n},\cdots,\, 2P_{m,n}\cos(m\lambda),\, 2P_{m,n}\sin(m\lambda) \cdots \Big) \textrm{  with  } m,n \leq N.
\end{equation*}
The TUKI starts with $\theta_0 \sim \N(\omega_b, Z_0 Z_0^T)$. For both approaches, we have $N_r = 63$. 
It is worth mentioning the inverse problem is ill-posed, another possible initial condition can be obtained by mirroring around the equator, since there is no observation in the south hemisphere. Hence, regularization is required, and we set the regularization parameter $\alpha$ to $0.5$, since $\omega_b$ is a good approximation.

\begin{figure}[ht]
\centering
\includegraphics[width=0.98\textwidth]{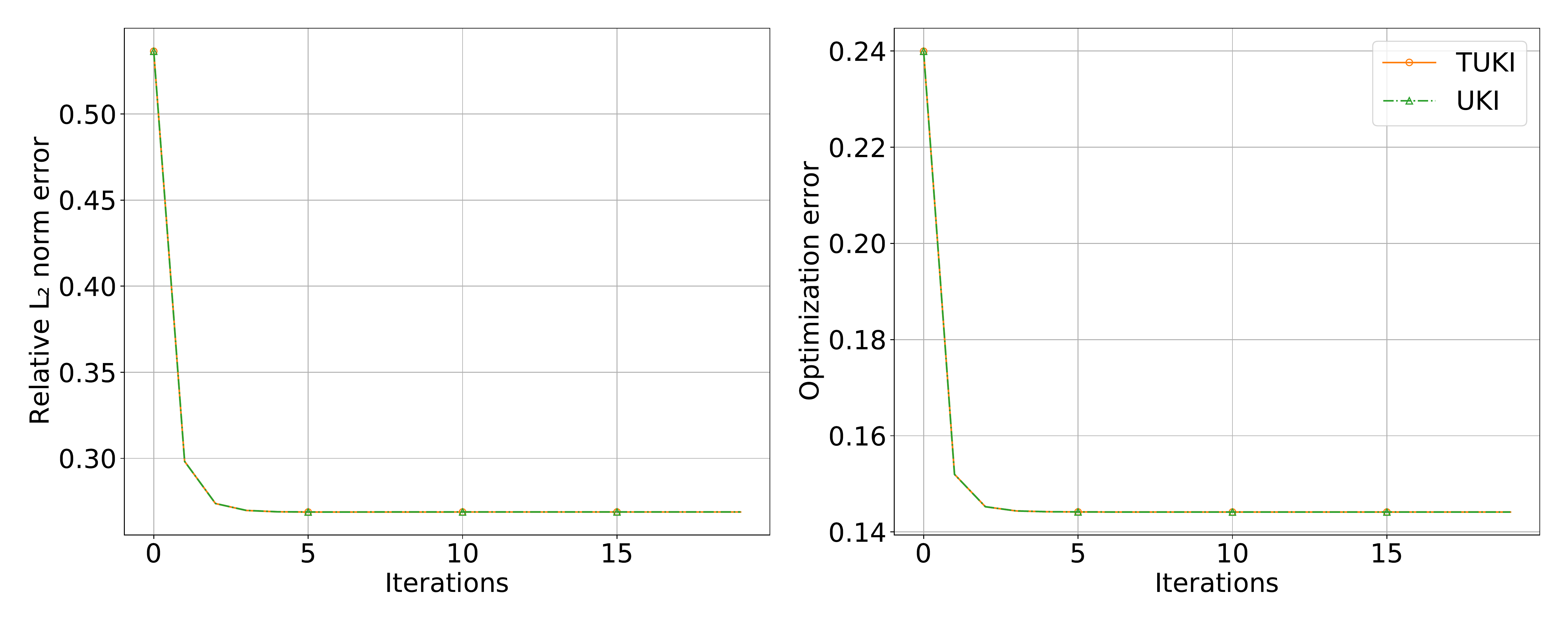}
\caption{$L_2$ error $\frac{\lVert\omega(x, \mean_n) - \omega_0\rVert_2}{\Vert\omega_0\rVert_2}$~(left) and the optimization error $\displaystyle \frac{1}{2}\lVert\Sigma_{\nu}^{-\frac{1}{2}} (y_{obs} - \py_n)\rVert^2$~(right) of the barotropic model problem.}
\label{fig:barotropic-converge}
\end{figure}

\begin{figure}[ht]
\centering
\includegraphics[width=0.48\textwidth]{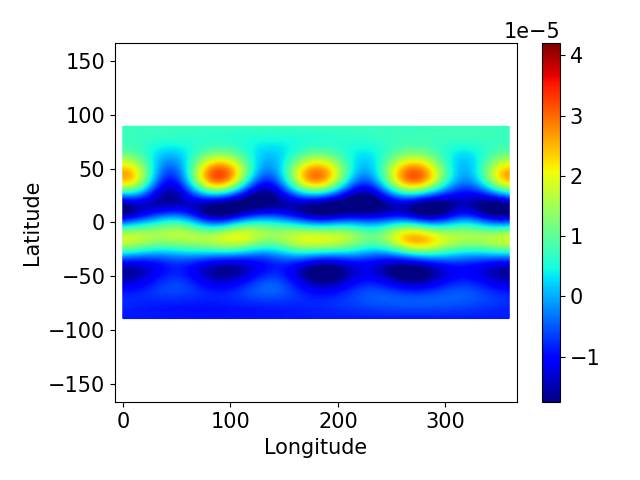}
\includegraphics[width=0.48\textwidth]{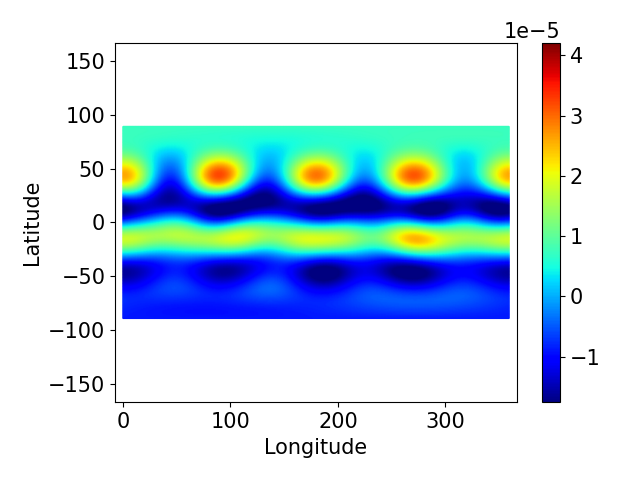}
\caption{Initial vorticity fields $\omega_0(x, \mean_n)$ recovered by UKI~(left) and TUKI~(right). The reference initial vorticity field is in~\cref{fig:vor_0}-d.}
\label{fig:barotropic-vorticity}
\end{figure}

\begin{figure}[ht]
\centering
\includegraphics[width=0.48\textwidth]{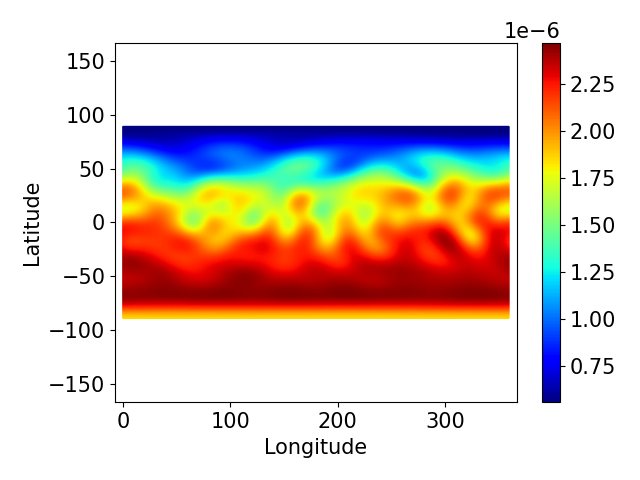}
\includegraphics[width=0.48\textwidth]{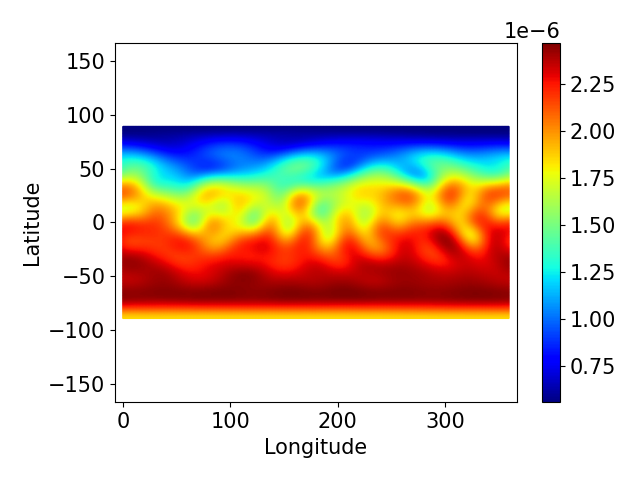}
\caption{Error variance fields in terms of the standard deviation at each point obtained by UKI~(left) and TUKI~(right).}
\label{fig:barotropic-vorticity-uq}
\end{figure}

The convergence of the initial vorticity field $\omega_0(x, \mean_n)$ and the optimization errors
at each iteration are depicted in \cref{fig:barotropic-converge}. 
TUKI performs very close to UKI with reparameterization for this case.  
The estimated initial vorticity fields $\omega_0(x,m_n)$ at the 20th iteration are depicted in \cref{fig:barotropic-vorticity}. 
 Both UKI and TUKI capture these sinusoidal disturbance of the truth initial field. The 3-$\sigma$ error fields at each point are depicted in \cref{fig:barotropic-vorticity-uq}. 
 Although the covariance  does
not represent the Bayesian posterior uncertainty, it does indicate the sensitivities inherent in the
estimation problem, and in particular  it features largest uncertainty in the south hemisphere, 
 since there is no observation.

\subsection{Idealized Generalized Circulation Model Problem}

\label{sec:app:GCM}
Finally, we consider using reduced-order model techniques to speed up an idealized general circulation model inverse problem. The model is based on the 3D Navier-Stokes equations, making the hydrostatic and shallow-atmosphere approximations common in atmospheric modeling. 
Specifically, we test on the notable Held-Suarez test case~\cite{held1994proposal}, in which a detailed radiative transfer model is replaced by Newtonian relaxation of temperatures toward a prescribed ``radiative equilibrium'' $T_{\mathrm{eq}}(\phi, p)$ that varies with latitude $\phi$ and pressure $p$. Specifically, the thermodynamic equation for temperature $T$ 
\[
\frac{\partial T}{\partial t} + \dots = Q
\]
(dots denoting advective and pressure work terms) contains a diabatic heat source 
\[
Q = -k_T(\phi, p, p_s) \bigl(T - T_{\mathrm{eq}}(\phi, p)\bigr),
\]
with relaxation coefficient (inverse relaxation time)
\[
k_T = k_a + (k_{s} - k_a)\max\Bigl(0, \frac{\sigma - \sigma_b}{1 - \sigma_b}\Bigr)\cos^4\phi.
\]
Here, $\sigma = p/p_s$, which is pressure $p$ normalized by surface pressure $p_s$, is the vertical coordinate of the model, and 
\[
T_{\mathrm{eq}} = \max\Bigl\{200K, \Bigl[315K - \Delta T_y \sin^2\phi - \Delta\theta_z\log\Bigl(\frac{p}{p_0}\Bigr)\cos^2\phi\Bigr]\Bigl(\frac{p}{p_0}\Bigr)^{\kappa}\Bigr\}
\]
is the equilibrium temperature profile ($p_0 = 10^5~\mathrm{Pa}$ is a reference surface pressure and $\kappa = 2/7$ is the adiabatic exponent). Default parameters are 
\[
k_a = (40\ \mathrm{day})^{-1},\qquad  
k_{s} = (4\ \mathrm{day})^{-1}, \qquad
\Delta T_y = 60\ \mathrm{K}, \qquad
\Delta\theta_z = 10\ \mathrm{K}.
\]

\begin{figure}[ht]
\centering
\includegraphics[width=0.48\textwidth]{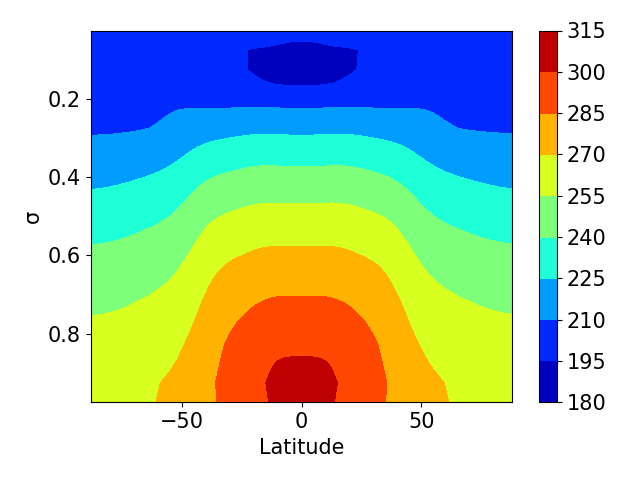}
\includegraphics[width=0.48\textwidth]{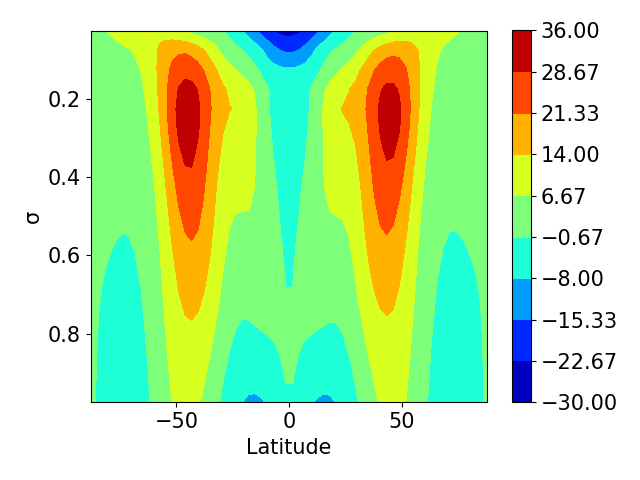}\\
\includegraphics[width=0.48\textwidth]{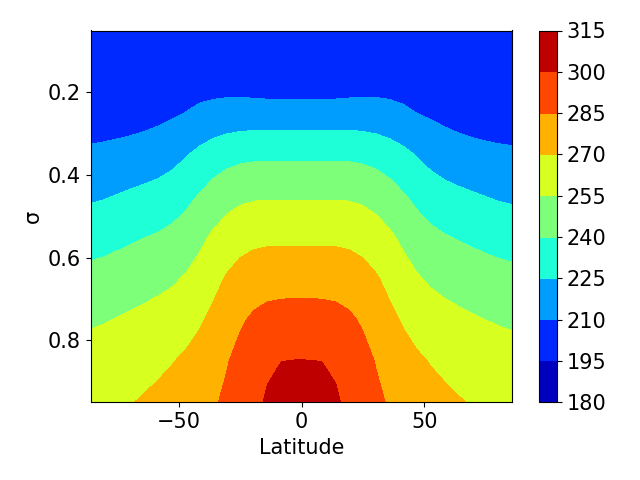}
\includegraphics[width=0.48\textwidth]{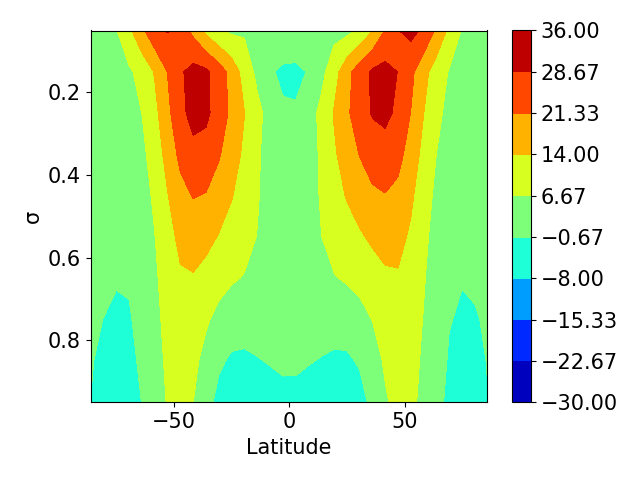}
\caption{Zonal mean temperature~(left) and zonal wind velocity~(right) obtained with the T42 grid~(top) and the T21 grid~(bottom).}
\label{fig:GCM-sol}
\end{figure}

For the numerical simulations, we use the spectral transform method in the horizontal, with T42 spectral resolution~(triangular truncation at wavenumber 42, with $64 \times 128$ points on the latitude-longitude transform grid); we use 20 vertical levels equally spaced in $\sigma$. With the default parameters, the model produce an Earth-like zonal-mean circulation, albeit without moisture or precipitation. A single jet is generated with maximum strength of roughly $30\ \mathrm{m~s^{-1}}$ near $45^{\circ}$ latitude (See \cref{fig:GCM-sol}).

Our inverse problem is constructed to learn parameters in Newtonian relaxation term $Q$:
\begin{equation*}
    (k_a,\ k_{s},\ \Delta T_y,\ \Delta\theta_z),
\end{equation*}
We do so in the presence of the following constraints: 
$$0\  \mathrm{day}^{-1}<k_a < 1\  \mathrm{day}^{-1}, \qquad   k_a<k_s <1\  \mathrm{day}^{-1}, \qquad   0\  \mathrm{K} < \Delta T_y, \qquad   0\  \mathrm{K}<\Delta\theta_z.$$
The inverse problem is formed as follows~\cite{UKI},
\begin{equation}
    y = \G(\theta) + \eta \quad \mathrm{with} \quad \G(\theta)=\overline{T}(\phi, \sigma)
\end{equation}
with the parameter transformation 
\begin{equation}
    \theta: (k_a, k_s, \Delta T_y, \Delta\theta_z) = \Big(\frac{1}{1+|\theta_{(1)}|},\ \frac{1}{1+|\theta_{(1)}|}+\frac{1}{1+|\theta_{(2)}|},\
    |\theta_{(3)}|,
    |\theta_{(4)}|\Big).
\end{equation}
The observation mapping $\G$ is defined by mapping from the unknown $\theta$ to 
the 200-day zonal mean of the temperature~($\overline{T}$) as a function of latitude ($\phi$) and height ($\sigma$), after an initial spin-up of 200 days.
The truth observation is the 1000-day zonal mean of the temperature~(see \cref{fig:GCM-sol}-top-left), after an initial spin-up 200 days to eliminate the influence of the initial condition. Because the truth observations come from an average 5 times as long as the observation window used for parameter learning, the chaotic internal variability of the model introduces noise in the observations. 

\begin{figure}[ht]
\centering
\includegraphics[width=0.7\textwidth]{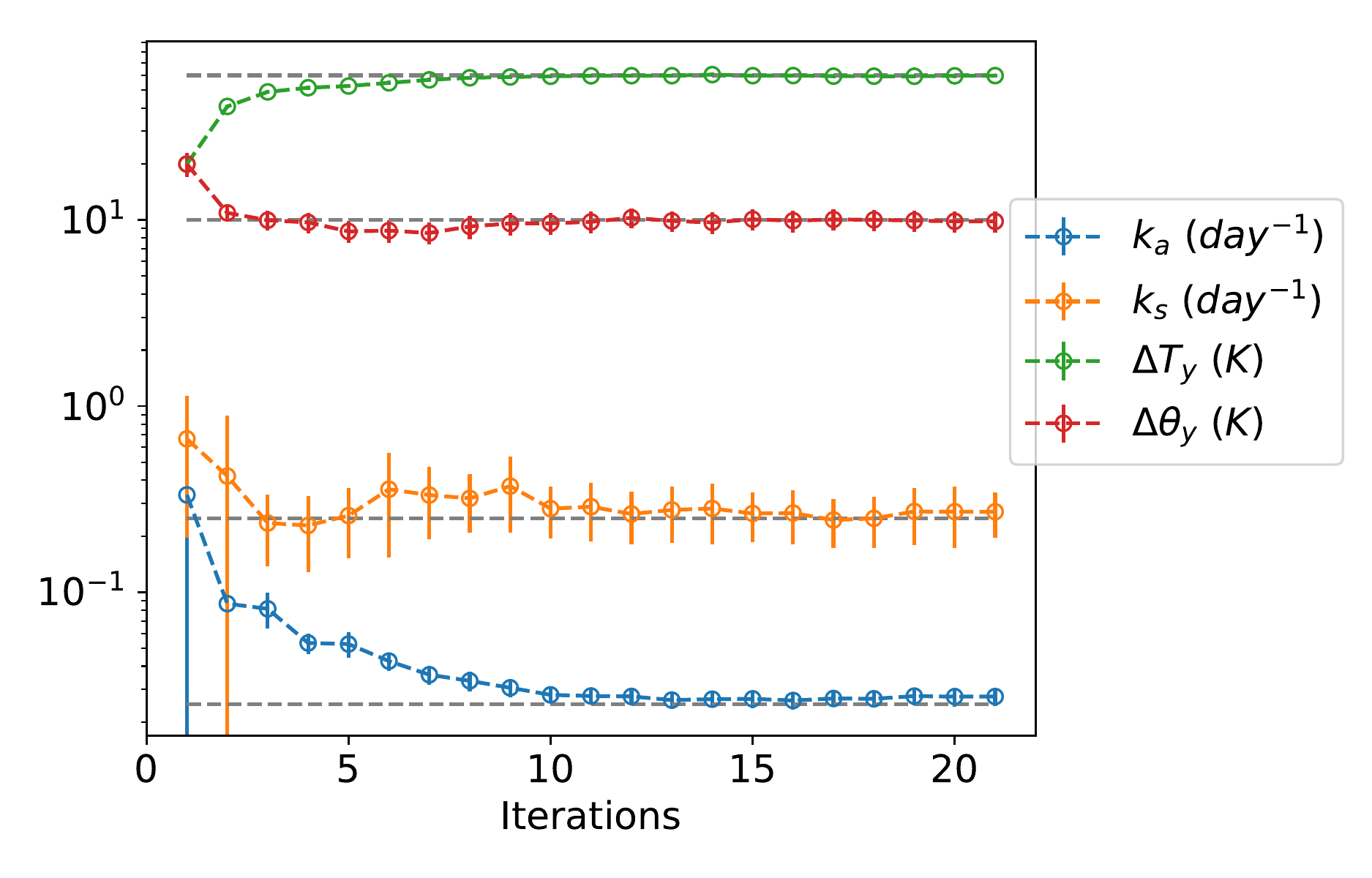}
\caption{Convergence of the idealized general circulation model inverse problem with UKI, true parameter values are represented by dashed grey lines.}
\label{fig:GCM-obj}
\end{figure}

To perform the inversion, we set $r_0 = [2\ \textrm{day},\ 2\ \textrm{day},\ 20\ \textrm{K},\ 20\ \textrm{K}]^{T}$  and UKI is initialized with $\displaystyle \theta_0 \sim \N\Big(r_0,\ \I\Big).$ Within the algorithm, we assume that the observation error satisfies $\eta \sim \N(0\ \textrm{K}, 3^2\I\ \textrm{K}^2)$. Because the problem is over-determined,
we set $\alpha=1.$ 
The reduced-order model technique discussed in~\cref{sec:speedup-rom} is applied to speed up the UKI. These $2N_{\theta}$  forward model evaluations are computed on a T21 grid~(triangular truncation at wavenumber 21, with $32 \times 64$ points on the latitude-longitude transform grid) with 10 vertical levels equally spaced in $\sigma$~(twice coarser in all three directions). The 1000-day zonal mean of the temperature and velocity predicted by the low-resolution model with the truth parameters are shown in~\cref{fig:GCM-sol}-bottom. Although there are large discrepancies comparing with results computed on the T42 grid, these low-resolution results only affect the covariance and $\F_u d\G$. 
The estimated parameters and associated $3-\sigma$ confidence intervals for each component at each iteration are depicted in \cref{fig:GCM-obj}. The estimation of model parameters at the 20th iteration are 
\begin{equation*}
\tiny
    \begin{bmatrix}
    k_a\\
    k_s\\
    \Delta T_y \\
    \Delta \theta_z
    \end{bmatrix}
    \sim
    \N\Big(
    \begin{bmatrix}
    0.0275\ \textrm{day}^{-1}\\
    0.271\ \textrm{day}^{-1}\\
    59.76\ \textrm{K}\\ 
    9.84\ \textrm{K}
    \end{bmatrix},\
    \begin{bmatrix}
    9.98\times10^{-7}\ \textrm{day}^{-2}  & -2.1\times10^{-6}\ \textrm{day}^{-2}  & -2.28\times10^{-4}\ \textrm{day}^{-1}\textrm{K} & 3.38\times10^{-5}\ \textrm{day}^{-1}\textrm{K}\\ 
    -2.1\times10^{-6}\ \textrm{day}^{-2}  & 6.05\times10^{-4}\ \textrm{day}^{-2}  & 5.22\times10^{-3}\ \textrm{day}^{-1}\textrm{K}  & -1.86\times10^{-3}\ \textrm{day}^{-1}\textrm{K}\\
    -2.28\times10^{-4}\ \textrm{day}^{-1}\textrm{K} & 5.22\times10^{-3}\ \textrm{day}^{-1}\textrm{K}  & 2.70\times10^{-1}\ \textrm{K}^{2}  & 3.19\times10^{-3}\ \textrm{K}^{2} \\ 
    3.38\times10^{-5}\ \textrm{day}^{-1}\textrm{K}  & -1.86\times10^{-3}\ \textrm{day}^{-1}\textrm{K} & 3.19\times10^{-3}\ \textrm{K}^{2}  & 1.85\times10^{-1}\ \textrm{K}^{2}
    \end{bmatrix}
    \Big).
\end{equation*}
The UKI with reduced-order models converges efficiently to the true parameters~(less than 20 iterations). Moreover,  the computational cost of these $2N_{\theta}$ low-resolution forward evaluations is even lower than a single high-resolution forward evaluation. 
Hence, for this inverse problems, the derivative-free UKI equipped with reduced order models achieves a CPU time speedup factor of $\bigO(10)$ and is as efficient as gradient-based optimization methods, which generally require an equally expensive backward propagation.

\section{Conclusion}
Unscented Kalman inversion is attractive for at least four main reasons: (i) it is derivative-free; (ii) it is robust for noise observations and chaotic inverse problems; (iii) it can be embarrassingly parallel; (iv) it provides sensitivity information.
 In the present work,  several strategies of making the UKI more efficient is presented for large scale inverse problems, including  using low-rank approximation and reduced-order model techniques.
Although we demonstrate the success of applying UKI for high dimensional inverse problems, the construction of low-rank basis for general scientific and engineering problems with complex geometries is not trivial, which is worth further investigation.  Another interesting area for future work would be to speed up UKI with other reduced-order model techniques.

\paragraph{Acknowledgments} D.Z.H. is supported by the generosity of Eric and Wendy Schmidt by recommendation of the Schmidt Futures program. J.H. is supported by the Simons Foundation as a Junior Fellow at New York University.

\appendix
\section{Proof of Theorems}
\label{sec:app:proof}
\begin{proof}[Proof of lemma~\ref{lemma:TUKI}]
We will prove that $\CovZ_n$ is rank $N_r$ and share the same column vector space spanned by the column vectors of $Z_0$ 
by induction.

When $n=0$, since 
$\textrm{rank}(\CovZ_0) = N_r$, therefore $\CovZ_n$ is rank $N_r$ and share the same column vector space spanned by the column vectors of $Z_0$.

We assume this holds for all $n \leq k$.
Since $\displaystyle Z_{\omega} = \sqrt{2-\alpha^2}Z_0$, the column space and the left singular vector space of $\Big(\alpha\CovZ_k \quad \CovZ_{\omega}\Big)$ are the same as the column space of $Z_0$.
Therefore the $N_r$-TSVD 
$$\Big(\alpha\CovZ_k \quad \CovZ_{\omega}\Big) = \widehat{U}_k \sqrt{\widehat{D}_k} \widehat{V}_k^T$$
is exact.
Next, since we have
\begin{equation*}
    \pCovZ_{k+1} \pCovZ_{k+1}^T = \Big(\alpha\CovZ_k \quad \CovZ_{\omega}\Big)\Big(\alpha\CovZ_k \quad \CovZ_{\omega}\Big)^T,
\end{equation*}
$\pCovZ_{k+1}$ is rank $N_r$ and the column space of $\pCovZ_{k+1}$ is the same as $Z_0$. 
Finally, since $\widehat{P}_{k+1}(\widehat{\Gamma}_{k+1}+\I)^{-1/2}$ is not singular, from 
$$\CovZ_{k+1} = \pCovZ_{k+1}\widehat{P}_{k+1}(\widehat{\Gamma}_{k+1}+\I)^{-1/2},$$
we have that 
$\CovZ_{k+1}$ shares the column space of  $\pCovZ_{k+1}$.

Therefore $\{\CovZ_n\}$ and $\{\pCovZ_n\}$ are rank $N_r$ and share the same column vector space spanned by the column vectors of $\CovZ_0$~(\ref{eq:TUKI-start}).
\end{proof}

\begin{proof}[Proof of theorem~\ref{theorem:negative}]
Gaussian vectors are rotation invariant. We can take any $N_\theta\times N_\theta$ orthogonal matrix $O$, then the joint law of $\{O\theta^1, O\theta^2,\cdots, O\theta^J\}$ is the same as the joint law of $\{\theta^1, \theta^2,\cdots, \theta^J\}$. Moreover, rotation preserves distance:
\begin{align}
    \textrm{dist}(\theta_{ref}, \textrm{span}\{\theta^1, \theta^2, \cdots, \theta^J\})=\textrm{dist}(O\theta_{ref}, \textrm{span}\{O\theta^1, O\theta^2, \cdots, O\theta^J\}),
\end{align}
which has the same law as $\textrm{dist}(O\theta_{ref}, \textrm{span}\{\theta^1, \theta^2, \cdots, \theta^J\})$. By taking expectation over $\theta^1, \theta^2, \cdots, \theta^J$ on both sides, we get
\begin{align}\label{e:eqO}
    \E[\textrm{dist}(\theta_{ref}, \textrm{span}\{\theta^1, \theta^2, \cdots, \theta^J\})^2]=\E[\textrm{dist}(O\theta_{ref}, \textrm{span}\{\theta^1, \theta^2, \cdots, \theta^J\})^2],
\end{align}
for any $N_\theta\times N_\theta$ orthogonal matrix $O$. We can further average \eqref{e:eqO}, over the uniform distribution on all $N_\theta\times N_\theta $ orthogonal matrices (Haar measure over the orthogonal group):
\begin{align}\label{e:EO}
 \E[\textrm{dist}(\theta_{ref}, \textrm{span}\{\theta^1, \theta^2, \cdots, \theta^J\})^2]=\E[\E_O[\textrm{dist}(O\theta_{ref}, \textrm{span}\{\theta^1, \theta^2, \cdots, \theta^J\})^2]],
\end{align}
where $\E_O$ is the expectation with respect to the uniform distribution on all $N_\theta\times N_\theta $ orthogonal matrices. Under this measure, $O\theta_{ref}$ is uniformly distributed on the sphere of radius $\|\theta_{ref}\|_2$, i.e. it has the same law as $\|\theta_{ref}\|_2\omega$, where $\omega$ is uniformly distributed on the unit sphere. We can rewrite the right hand side of \eqref{e:EO} as
\begin{align}\label{e:outO}
  \E_O[\textrm{dist}(O\theta_{ref}, \textrm{span}\{\theta^1, \theta^2, \cdots, \theta^J\})^2]
  =\|\theta_{ref}\|^2_2\E_\omega[\textrm{dist}(\omega, \textrm{span}\{\theta^1, \theta^2, \cdots, \theta^J\})^2].
\end{align}
For Gaussian vectors $\theta^1, \theta^2,\cdots, \theta^J$, almost surely, their span, $\text{span}\{\theta^1, \theta^2, \cdots, \theta^J\}$ has dimension $J$.  $\textrm{dist}(\omega, \textrm{span}\{\theta^1, \theta^2, \cdots, \theta^J\})$ is the distance from 
a uniformly distributed random vector on the unit sphere to a plane of dimension $J$, we can again rotate the plane, which will not change the law. So we can rotate the plane to the plane spanned by the coordinate vectors $e_1, e_2, \cdots, e_J$, then 
\begin{align}\begin{split}\label{e:Ew}
    \E_\omega[\textrm{dist}(\omega, \textrm{span}\{\theta^1, \theta^2, \cdots, \theta^J\})^2]
    &=\E_\omega[\textrm{dist}(\omega, \textrm{span}\{e_1, e_2, \cdots, e_J\})^2]\\
    &=\E_\omega\left[\sum_{i=J}^{N_\theta}\omega_i^2\right]=\frac{N_\theta-J}{N_\theta}=1-\frac{J}{N_\theta},
\end{split}\end{align}
where we use that for a uniformly distributed vector, the expectation of each coordinate square is $1/N_\theta$. The claim \eqref{e:lob} follows from combining \eqref{e:outO} and \eqref{e:Ew}.



\end{proof}

\section{Ensemble Based Kalman Inversion}
\label{sec:app:ensemble_method}
The stochastic ensemble Kalman inversion~\cite{UKI} is 
\begin{itemize}
  \item Prediction step : \begin{align*}
                      &\pp_{n+1}^{j} = \alpha \theta_{n}^{j} + (1 - \alpha)r  + \omega_{n+1}^{j} \qquad \pmean_{n+1} = \frac{1}{J}\sum_{j=1}^{J}\pp_{n+1}^{j}, \\
                       &\pCov_{n+1} = \frac{1}{J-1} \sum_{j=1}^{J} (\pp_{n+1}^{j} - \pmean_{n+1})(\pp_{n+1}^{j}-\pmean_{n+1} )^T.\\
                   \end{align*}
  \item  Analysis step : 
  \begin{equation}
  \label{eq:EKI-analysis}
  \begin{split}
         &\ppy_{n+1}^{j} = \G(\pp_{n+1}^{j})  \qquad \py_{n+1} = \frac{1}{J}\sum_{j=1}^{J}\ppy_{n+1}^{j},\\
        &\pCov_{n+1}^{\theta p} = \frac{1}{J-1}\sum_{j=1}^{J}(\pp_{n+1}^{j} - \pmean_{n+1})(\ppy_{n+1}^{j} - \py_{n+1})^T,  \\
        &\pCov_{n+1}^{p p} = \frac{1}{J-1}\sum_{j=1}^{J}(\ppy_{n+1}^{j} - \py_{n+1})(\ppy_{n+1}^{j} - \py_{n+1})^T +\Sigma_{\nu}, \\
        &\theta_{n+1}^{j} = \pp_{n+1}^{j} + \pCov_{n+1}^{\theta p}\left(\pCov_{n+1}^{pp}\right)^{-1}(y - \ppy_{n+1}^{j} - \nu_{n+1}^{j}),\\
        &\mean_{n+1} = \frac{1}{J} \sum_{j=1}^{J} \theta_{n+1}^{j}.\\
  \end{split}
\end{equation}
\end{itemize}
Here the superscript $j = 1,\cdots,\ J$ is the ensemble particle index, $\omega_{n+1}^{j} \sim \N(0, \Sigma_{\omega})$ and  $\nu_{n+1}^{j} \sim \N(0, \Sigma_{\nu})$ are independent and identically distributed random variables.

Inspired by square root Kalman filters~\cite{bishop2001adaptive, anderson2001ensemble,tippett2003ensemble, wang2003comparison},  the analysis step of the stochastic Kalman inversion can be replaced by a deterministic analysis, which leads to deterministic ensemble Kalman inversions.
Deterministic ensemble Kalman inversions update the mean first, as following,
\begin{equation}
     \mean_{n+1} = \pmean_{n+1} + \pCov_{n+1}^{\theta p}\left(\pCov_{n+1}^{pp}\right)^{-1}(y - \py_{n+1}).\\
\end{equation}
Then we can define the matrix square roots $\pCovZ_{n+1},\, \CovZ_{n+1} \in \R^{N_{\theta}\times J}$ of $\pCov_{n+1},\,\Cov_{n+1}$ as following,
\begin{equation}
\begin{split}
    \pCovZ_{n+1} = \frac{1}{\sqrt{J-1}}\Big(\pp_{n+1}^{1} - \pmean_{n+1}\quad \pp_{n+1}^{2} - \pmean_{n+1}\quad...\quad\pp_{n+1}^{J} - \pmean_{n+1} \Big),\\
     \CovZ_{n+1} = \frac{1}{\sqrt{J-1}}\Big(\theta_{n+1}^{1} - \mean_{n+1}\quad \theta_{n+1}^{2} - \mean_{n+1}\quad...\quad\theta_{n+1}^{J} - \mean_{n+1} \Big).
     \end{split}
\end{equation}
The particles $\{\theta_{n+1}^{j}\}$ can be deterministically updated by the ensemble adjustment/transform Kalman methods, which brings about following ensemble adjustment/transform Kalman inversions.

\subsection{Ensemble Adjustment Kalman Inversion}
Following the ensemble adjustment Kalman filter proposed in~\cite{anderson2001ensemble}, the analysis step updates particles deterministically with a pre-multiplier $A$,
\begin{equation}
\theta_{n+1}^j - \mean_{n+1} = A (\pp_{n+1}^{j} - \pmean_{n+1}). \\
\end{equation}
Here $A = F \sqrt{D^{p}} U \sqrt{D}\sqrt{D^{p}}^{-1}F^T $ with 
\begin{equation}
\label{eq:eaki-svd}
\begin{split}
   \textrm{SVD :}     \quad       &\pCovZ_{n+1} =  F \sqrt{D^p} V^T,\\
   \textrm{SVD :}     \quad      &V^T\Big(\I + \pGCovZ_{n+1}^T \Sigma_{\nu}^{-1}  \pGCovZ_{n+1}\Big)^{-1} V = U D U^T,
\end{split}
\end{equation}
where both $\sqrt{D^p}$ and $D$ are non-singular diagonal matrices, with dimensionality rank($\pCovZ_{n+1}$) and 
\begin{equation}
\label{eq:EAKI-V}
    \pGCovZ_{n+1} = \frac{1}{\sqrt{J-1}}\Big(\py_{n+1}^{1} - \py_{n+1}\quad \py_{n+1}^{2} - \py_{n+1}\quad...\quad\py_{n+1}^{J} - \py_{n+1} \Big).
\end{equation}
It can be verified that the ensemble covariance matrix from~\cref{eq:EKI-analysis} satisfies
\begin{equation}
\label{eq:cov-eaki}
    \begin{split}
        \Cov_{n+1} &= \CovZ_{n+1} \CovZ_{n+1}^T \\ 
               &= \pCovZ_{n+1}\Big(\I -  \pGCovZ_{n+1}^T (\pGCovZ_{n+1}\pGCovZ_{n+1}^T + \Sigma_{\nu})^{-1}\pGCovZ_{n+1}\Big)\pCovZ_{n+1}^T \\
               &= \pCovZ_{n+1}\Big(\I +  \pGCovZ_{n+1}^T  \Sigma_{\nu}^{-1}\pGCovZ_{n+1}\Big)^{-1}\pCovZ_{n+1}^T \\
    \end{split}
\end{equation}
Bringing~\cref{eq:eaki-svd} into~\cref{eq:cov-eaki} leads to 
\begin{equation}
\label{eq:cov-eaki-2}
    \begin{split}
        \Cov_{n+1} 
               &= F \sqrt{D^{p}}  U D U^{T} \sqrt{D^p}F  \\
               & = F \sqrt{D^{p}} U \sqrt{D} \sqrt{D^p}^{-1}F^T \pCov_{n+1} F \sqrt{D^p}^{-1} \sqrt{D} U^{T}  \sqrt{D^p}F^T\\
               & = A\pCovZ_{n+1}  \pCovZ_{n+1}^T A^T.
    \end{split}
\end{equation}
And therefore the ensemble adjustment Kalman filter delivers the same covariance matrix as the stochastic ensemble Kalman filter. Although the covariance matrix is not required in both algorithms. 

\subsection{Ensemble Transform Kalman Inversion}
Following the ensemble transform Kalman filter proposed in~\cite{bishop2001adaptive}, the analysis step updates particles deterministically with a post-multiplier $T$,
\begin{equation}
\CovZ_{n+1} = \pCovZ_{n+1} T.\\
\end{equation}
Here $T = P(\Gamma + I)^{-1/2}$, with 
\begin{equation}
\textrm{SVD:} \quad \pGCovZ_{n+1} \Sigma_{\nu}^{-1} \pGCovZ_{n+1} = P\Gamma P^T,
\end{equation}
where $\pGCovZ_{n+1}$ is defined in~\cref{eq:EAKI-V}.
Following \cref{eq:cov-eaki}, the ensemble covariance matrix satisfies
\begin{equation}
\begin{split}
\Cov_{n+1} &= \pCovZ_{n+1}\Big(\I +  \pGCovZ_{n+1}^T  \Sigma_{\nu}^{-1}\pGCovZ_{n+1}\Big)^{-1}\pCovZ_{n+1}^T \\
           &= \pCovZ_{n+1}\Big(\I +  P\Gamma P^T\Big)^{-1}\pCovZ_{n+1}^T \\
           &= \pCovZ_{n+1} P\Big(\I +  \Gamma \Big)^{-1}P^T\pCovZ_{n+1}^T \\
           &= \pCovZ_{n+1} P\Big(\I +  \Gamma \Big)^{-1}P^T\pCovZ_{n+1}^T \\
           & = \pCovZ_{n+1} T T^T \pCovZ_{n+1}^T.
\end{split}
\end{equation}  
And therefore the ensemble transform Kalman filter delivers the same covariance matrix as the stochastic ensemble Kalman filter. Although the covariance matrix is not required in both algorithms. 

\begin{remark}
The original ETKF is biased, since $\CovZ_{n+1} 1 \neq 0$~($\mean_{n+1}$ is not the mean of $\{\theta_{n+1}^{j}\}_j$). 
An unbiased ETKF fix is introduced in \cite{wang2003comparison} by defining a symmetric post-multiplier $T = P(\Gamma + \I)^{-1/2}P^T$. Since 
\begin{equation}
    TT^T 1= \Big(\I -  \pGCovZ_{n+1}^T (\pGCovZ_{n+1}\pGCovZ_{n+1}^T + \Sigma_{\nu})^{-1}\pGCovZ_{n+1}\Big) 1 = 1,
\end{equation}
$1$ is the eigenvector of $T$ and, therefore 
\begin{equation}
    \CovZ_{n+1} 1 = \pCovZ_{n+1} T 1  = \pCovZ_{n+1}1 = 0.
\end{equation}
\end{remark}


\bibliographystyle{unsrt}
\bibliography{references}
\end{document}